\documentclass[10pt,a4paper]{article}
\usepackage[utf8]{inputenc}
\usepackage{url,amsmath,amssymb,amsthm,amsfonts,epsfig,eurosym,geometry}
\usepackage{cite}

\def\Pr{\mathbb{P}}
\def\Es{\mathbb{E}}

\def\bq{\begin{quote}\begin{em}}
\def\eq{\end{em}\end{quote}}

\def\un{\mathbf{1}}

\newcommand{\R}{\mathbb R}
\newcommand{\N}{\mathbb N}
\newcommand{\Z}{\mathbb Z}

\newtheorem{theo}{Theorem}[section]
\newtheorem{coro}{Corollary}[section]
\newtheorem{lemm}{Lemma}[section]
\newtheorem{defi}[theo]{Definition}
\newtheorem{prop}[theo]{Proposition}
\newtheorem{rque}[theo]{Remark}
\newtheorem{exem}[theo]{Example}
\newtheorem{exer}[section]{Exercice}

\def\demi{\frac{1}{2}}
\def\ds{\displaystyle}
\def\btheo{\begin{theo}}\def\etheo{\end{theo}}
\def\bcoro{\begin{coro}}\def\ecoro{\end{coro}}
\def\blemm{\begin{lemm}}\def\elemm{\end{lemm}}
\def\bdefi{\begin{defi}}\def\edefi{\end{defi}}
\def\bprop{\begin{prop}}\def\eprop{\end{prop}}
\def\brque{\begin{rque}}\def\erque{\end{rque}}
\def\bexem{\begin{exem}}\def\eexem{\end{exem}}
\def\bexer{\begin{exer}\begin{em}}\def\eexer{\end{em}\end{exer}}

\usepackage{graphics}
\usepackage{color}              
\usepackage{epsfig}

\def\1{{\rm 1\kern-.8ex 1}}

\begin{document}
\title{Optimal stopping time and halting set for total variation distance}
\author{}
\maketitle

\begin{abstract} 
An aperiodic and irreducible  Markov chain on a finite state space converges to its stationary distribution. When  convergence to equilibrium is measured by total variation distance, there exists an optimal coupling and a maximal coupling time. In this article, the maximal coupling time is compared to the hitting time of a specific state or set. Such sets, named halting sets, are studied in the case of symmetric birth-and-death chains and in some other examples. Some applications to the cutoff phenomenon are given. These results yield new methods to calculate cutoff times for some monotone birth-and death chains without the lazy hypothesis . 
\end{abstract}


\section{Introduction}
Let $\mathbb{X}$ be a finite space, $(X_t)_{t\geq 0}$ a discrete or continuous time irreducible Markov chain on $\mathbb{X}$. Let $\pi$ be its stationary distribution . For $t\geq 0$, denote by $\pi_t$ the distribution of $X_t$. Let $d$ be a distance between probability distributions on $\mathbb{X}$. For $t\geq 0$, denote by $d(t)$ the distance between the distribution of $X_t$ and $\pi$: 
\[d(t)=d(\pi_t,\pi)\]
For $A\subset \mathbb{X}$, denote by $T_A$ the hitting time of $A$.

\bdefi\label{def:haltingset}
A set $A\subset \mathbb{X}$ is a halting set for the distance $d$ and the Markov chain $(X_t)_{t\geq 0}$ if for all $t\geq 0$,
\[d(t)\leq\Pr(T_A>t)\].
\edefi

If the Markov chain is also aperiodic in the discrete time case, $d(t)$ goes to $0$ when $t$ goes to infinity. The question of non-asymptotic behavior is hard. There are several techniques for obtaining explicit bounds on $d(t)$ : Fourier analysis, coupling, strong stationnary times, $\ldots$

Several of them use hitting times.

Strong stationary times are allied with separation distance, $s(t)=\max_{x\in X}(1- \frac{\pi_t(x)}{\pi(x)})$. $T$ a randomized stopping time  is a strong stationary time if $T$ and $X_T$ are independent and if the law of $X_T$ is given by the stationary distribution $\pi$. In  \cite{AldousDiaconis87}, Aldous and Diaconis prove that there exists a strong stationary time $T$ such that for all $t\geq 0$, $s(t)=\Pr(T>t)$. Furthermore, in \cite{DiaconisFill90}, Diaconis and Fill construct an absorbing dual Markov chain $X^*$ such that $T$ is the hitting time of an absorbing state for $X^*$. This is used by Diaconis and Saloff-Coste in \cite{DiaconisSaloff2} for birth-and-death chain to prove some cut-off results.

Let $T$ be a strong stationary time. A state $y^*$ is a halting state for $T$ if $X_t=y^*$ implies $T\leq t$ or equivalently if $T\leq T_{y^*}$. Furthermore, if there exists a halting state for $T$, then  $s(t)=\Pr(T>t)$ for all $t$. This is an example of the use of hitting times to obtain bounds on a distance.

In coupling methods, hitting times are very often used to bound $\Pr(T>t)$ where $T$ is the coupling time. Several examples can be found in the book of Levin, Peres and Wilmer \cite{LevinPeresWilmer} p.65.

In \cite{MartinezYcart}, Martinez and Ycart also use hitting times to obtain bounds on the total variation distance. They prove that in general, the access time to equilibrium and the hitting times tend to be equivalent, if the process starts ``far away''. But the context is different, they consider continuous time Markov chain on a countable set $I$ and they study $\|\delta_aP_t-\pi\|_{TV}$ when $a$ goes to infinity.

In this article, mainly the total variation distance is used, that is:
\[d(t)=\|\pi_t-\pi\|_{TV}:=\ds\sup_{A\subset\mathbb{X}}|\pi_t(A)-\pi(A)|=\demi\ds\sum_{x\in \mathbb{X}}|\pi_t(x)-\pi(x)|\]

Find some bounds on $d(t)$ is equivalent to bound the mixing time of the chain. For $\epsilon \in ]0,1[$, denote by $t_{\mbox{\it mix}}(\epsilon)$, defined by  
\[t_{\mbox{\it mix}}(\epsilon)=\min\{t\geq 0 / d(t)\leq \epsilon\}.\]

In \cite{ LevinPeresWilmer} (theorem 10.14 p 134), the autors established the following theorem for a lazy chain :
\[ t_{\mbox{\it mix}}(\frac{1}{4})\leq 2\max_{x\in X}\Es_{\pi}[T_x]+1\]
In fact, they prove that $d(t)\leq \sqrt{\frac{\Es_{\pi}[T_x]}{8t}}$ where $\pi_0=\delta_x$.

In \cite{DingLubetzkyPeres}, Ding, Lubetzky and Peres prove the following result for a lazy birth-and-death chain on $\{0,\ldots, N\}$ started at $0$,
\begin{equation}\label{eq:inequalityDLP} 
\forall 0\leq l\leq N,\ \forall n\in\N,\ \|\pi_n-\pi\|_{TV}\leq \Pr(T_l>n)+\pi(\{l+1,\ldots,N\})
\end{equation}

So it is a natural idea and often used method to compare the access time to equilibrium with the hitting time of a given state. This state will depend on the initial state.

One object of this article is to look at ``minimal'' coupling time for the total variation distance, mimicking in some meaning what is done for the separation distance. The construction of Aldous and Diaconis of a randomized stopping time satisfying $s(n)=\Pr(T>n)$ for all $n$ leads naturally to a strong stationary stopping time. So rewriting what can be found in \cite{Griffeath1975} and \cite{Lindvall}, a randomized stopping time $T$ is constructed satisfying $d(t)=\Pr(T>t)$ for all $t$ and the properties of the law of $(X_T,T)$ lead to the following result:

\btheo\label{theo:condsufhalting}
Let ${\cal M}=\{x\in X /\  \forall t\geq 0,\ \pi_t(x)\leq \pi(x)\}$, then ${\cal M}$ is a halting set for the total variation distance and  the  chain $(X_t)_{t\geq 0}$.
\etheo

In the following, we  omit to say {\it for the total variation distance and  the  chain $(X_t)_{t\geq 0}$}.
Of course if $x^*\in{\cal M}$, $\{x^*\}$ is a halting set and so for all $t\geq 0$, $d(t)\leq \Pr(T_{x^*}>t)$, a property to be compared to (\ref{eq:inequalityDLP}). A such state is called a halting state.

The difficulty then will be to find halting sets. The first fact is  that if we have some property of monotonicity, a halting set must exist (proposition (\ref{prop:monotony}).

Halting sets can be used to prove some results about cutoff. This concepts was intoduced by Aldous and Diaconis in  \cite{AldousDiaconis86} to describe the fact that many ergodic Markov chains converge abruptly to their stationary distribution.

Consider a family $(X^{(N)}_n)_{n\geq 0}$ of aperiodic irreducible Markov chains on a finite state space $X^{(N)}$, each with its stationary distribution $\pi^{(N)}$and its distance from stationary $d^{(N)}_n$. $t_N$ is a cutoff if
\[\begin{array}{ll}
 d^{(N)}([ct_N]) \longrightarrow 0 \ \ &\mbox{if }\ \ c>1\\
d^{(N)}([ct_N]) \longrightarrow 1 \ \ &\mbox{if }\ \ c<1
  \end{array}\]
  
Diaconis and Saloff-Coste in\cite{DiaconisSaloff2} verified this conjecture for continuous-time birth-and-death chains, started at an endpoint, with convergence measured in separation.

Ding, Lubetzky and Peres in \cite{DingLubetzkyPeres} proved it for continuous-time birth-and-death chains and lazy discrete-time birth-and-death chains, started at an endpoint, with convergence measured in total variation. For such chains, the fact that the product of the mixing time and the spectral gap tends to infinity is equivalent to the fact that the product of the expected value of the hitting time of the median of the stationary distribution and the spectral gap tends to infinity.

The cutoff depends of the distance. For example, the  Ehrenfest process starting at $0$ on $\{0,\ldots, N\}$ has a $\frac{N\ln N}{2}$ separation cutoff but a $\frac{N\ln N}{4}$ total variation and $L^2$ cutoff. The biased (p; q)-random walk starting at $0$ on $\{0,\ldots, N\}$ has a $\frac{N}{p-q}$ total variation and separation cutoff but a $\frac{\ln(\frac{p}{q})}{2(1-2\sqrt{pq})}N$ $L^2$ cutoff.

Recently, Basu, Hermon and Peres in \cite{BasuHermonPeres}, study the link between cutoff and some concentration of hitting time of ``worst'' sets of stationary measure at least $\alpha$, for some $\alpha\in ]0,1[$. They prove that in the case of a lazy reversible irreducible Markov chain, cutoff is equivalent to a notion of cutoff for hitting times, denoted $hit_\alpha$-cutoff. Here, we use explicit hitting time, the hitting time of a halting state for total variation.

Recall the definition of a cutoff window :

If  $(w_N)$ and $(t_N)$ are two sequences such that $w_N=o(t_N)$, one may define that the family of chain $(X^{(N)}_n)_{n\geq 0}$ exhibits a cutoff at $t_N$ with window $w_N$ if
\[\begin{array}{l}
   \lim_{\gamma\rightarrow+\infty}\varliminf_{N\rightarrow+\infty}d^{(N)}([t_N-\gamma w_N])=1\\
\lim_{\gamma\rightarrow+\infty}\varlimsup_{N\rightarrow+\infty}d^{(N)}([t_N+\gamma w_N])=0

  \end{array}\]

A general result can be the following in the discrete time:
\bprop\label{prop:cutoffgeneral}
Let $A^{(N)}$ be a halting set such that $\sigma(T_{A^{(N)}})=o(\Es[T_{A^{(N)}}])$ . Let $\lambda_{A^{(N)}}$ be the largest eigenvalue of the restriction of the transition kernel of the chain to $\mathbb{X}\setminus A^{(N)}$.

If $(1-\lambda_{A^{(N)}})\Es[T_{A^{(N)}}]$ goes to infinity when $N$ goes to infinity, then the chain has a  cutoff at time $\Es[T_{A^{(N)}}]$ for total variation.

If furthermore, $\varliminf(1-\lambda_{A^{(N)}})\sigma(T_{A^{(N)}})>0$, then the chain has a  cutoff at time $\Es[T_{A^{(N)}}]$ with window $\sigma(T_{A^{(N)}})$ for total variation.
\eprop

The rest of the paper is to determine  halting set in several examples, essentialy in the case of birth-and-death chains. A Markov  chain on  $\{0,\ldots,2N+1\}$ on $\{0,\ldots,N\}$ with transitions given by $K=(k(x,y))_{x,y\in X}$ is symmetric if for all $x,y\in X$, $k(N-x,N-y)=k(x,y)$. An example of result is the following: In the case of a monotone  and symmetric birth-and-death chain on  $\{0,\ldots,2N+1\}$, $N+1$ is a halting state for total variation. Furthermore, if $1=\lambda_0>\lambda_1^{(N)}>\cdots>\lambda_{2N+1}^{(N)}$ are the eigenvalues of the transition kernel and if $(1-\lambda_1^{(N)})\displaystyle\sum_{k=1}^{2N+1}\frac{1}{1-\lambda_k^{(N)}}$ goes to infinity when $N$ goes to infinity, then the chain has a  cutoff at time $\demi \displaystyle\sum_{k=1}^{2N+1}\frac{1}{1-\lambda_k^{(N)}}$ for total variation.

The paper is organized as follows. In section 2, the construction of an  optimal randomized stopping time for the total variation distance for  discrete and continuous times markov chains is given. Section 3 exposes the notion of a halting set for total variation for a chain with an initial distribution and gives applications to cutoff. In section 4, general results about monotone birth-and-death Markov chain are given and the case of symmetric monotone birth-and-death Markov chain is studied. These results are used to prove cutoff in several examples. 

\section{Optimal  stopping time}\label{sec:optstopti}
\subsection{Construction of an optimal stopping time in discrete time}

This construction is not new but the presentation made here is different from Lindvall's proof in \cite{Lindvall}.

Let $(X_n)_{n\geq 0}$ be a Markov chain on a finite space $\mathbb X$ with a stationary probability $\pi$ such that for all $x\in X$, $\pi(x)>0$.

The transition kernel is denoted by $K=(k(x,y))_{x,y\in X}$. Let $\pi_0$ be the distribution of $X_0$ and for $n\geq 1$, $\pi_n=\pi_0K^n$ be the distribution of  $X_n$.
Let $d_n=\|\pi_n-\pi\|_{TV}$ be the total variation distance between $\pi_n$ and $\pi$.

We define for $n\geq 0$, $\gamma_n :X\rightarrow \R_+$ by
\[\gamma_0(x)=\pi_0\wedge \pi(x) \ \ \ \mbox{and for } n\geq 1,\ \gamma_n(x)=(\pi_n\wedge \pi)(x)-(\pi_{n-1}\wedge \pi)K(x)\]

As $d(n)=1-\ds\sum_{x\in X}(\pi_n(x)\wedge\pi(x))$, for all $n\geq 1$, $d(n-1)-d(n)=\ds\sum_{x\in X}\gamma_n(x)$.

By consequence, we are looking for a randomized stopping time $T$ which satisfies for all $n\geq 0$, for all $x\in X$, $\Pr(X_n=x,T=n)=\gamma_n(x)$.

\bprop \label{prop:defpsi}
Let $\Delta_n$  be defined by 
\begin{equation}\label{delta}\left\lbrace \begin{array}{l} \Delta_0=\pi_0\\
                     \Delta_{n+1}=\pi_{n+1}-(\pi_{n}\wedge\pi)K\
                    \end{array}\right. \end{equation}

Let $\Psi_n$  be defined for all $n\geq 0$ by 
\begin{equation}\label{psi}\psi_n(x)=\left\lbrace\begin{array}{cl}
                                                  \frac{\gamma_n(x)}{\Delta_n(x)}&\mbox{ if }\Delta_n(x)\not=0\\
                                                   1&\mbox{ if }\Delta_n(x)=0
                                                 \end{array}
\right.
\end{equation}

For all $n\geq 0$, $\psi_n$ takes value in $[0,1]$.

Let  $(U_n)_{n\geq 0}$ be independent variables uniformly distributed on $[0,1]$, independent of $(X_n)_{n\geq 0}$. 

Let $T$ be the randomized stopping time defined by $T=\inf\{n\geq 0\ /\ U_n\leq \psi_n(X_n)\}$.

 Then $\forall n\geq 1,\ \forall x\in {\mathbb X}$,
\begin{equation}\label{loiXTT}
 \Pr(X_n=x,T=n)=\gamma_n(x)
\end{equation} 
\begin{equation}\label{loiXTT2}
 \Pr(X_n=x,T\geq n)=\Delta_n(x), \Pr(X_n=x,T>n)=\pi_n(x)-(\pi_n\wedge \pi)(x)
\end{equation} 
\begin{equation}\label{loiT}
 \Pr(T>n)=d(n)
\end{equation} 

\eprop

\begin{proof}
As  $0\leq \gamma_n\leq \Delta_n$, $\psi_n$ takes value in $[0,1]$ and $T$ is well defined.

We prove (\ref{loiXTT2}) inductively :
\[\Delta_0(x)=\pi_0(x)=\Pr(X_0=x,T\geq 0)\]
\[\begin{array}{ll}
   \Pr(X_{n+1}=x,T\geq n+1)&=\sum_{y\in {\mathbb X}}\Pr(X_n=y,X_{n+1}=x,T\geq n,U_n>\psi_n(y))\\
&=\sum_{y\in {\mathbb X}}(1-\psi_n(y))\Pr(X_n=y,T\geq n)k(y,x)\\
&=((1-\psi_n)\Delta_n)K(x)=(\Delta_n-\psi_n\Delta_n)K(x)\\
&= (\Delta_n-\gamma_n)K(x)
=\Delta_{n+1}(x)
  \end{array}
\]
The second equality comes from the Markov property and definition of $T$.
(\ref{loiXTT}) and (\ref{loiT}) are easy consequences of (\ref{loiXTT2}).

\end{proof}

\brque
$\psi_n$ is given by the following formulae:
\begin{equation} \label{eq:psin}
 \psi_n(x)=\left\lbrace\begin{array}{cl}\frac{\gamma_n(x)}{\Delta_n(x)}=\frac{\pi_n\wedge\pi(x)-(\pi_{n-1}\wedge\pi)K(x)}{\pi_n(x)-(\pi_{n-1}\wedge\pi)K(x)}&\mbox{if }\Delta_n(x)\not=0\\
1&\mbox{if }\Delta_n(x)=0
                      \end{array}\right.
\end{equation} 
\erque

\brque
An alternative definition can be made for the randomized stopping time $T$.

We define for all $n\geq 0$
\begin{equation} \label{Jn}
 J_n=\prod_{k=0}^{n}(1-\psi_k(X_k))
\end{equation} 
and $T=\inf\{n\geq 0/\ U\geq J_n\}$ where $U$ is a random variable uniformly distributed $[0,1]$ independent of $(X_n)_{n\geq 0}$.
\erque

\brque
To study distance from stationarity for Markov chains, separation distance has good properties.

The separation distance is defined by 

\[ s(n):=sep(\pi_n,\pi)=\sup_{y\in {\mathbb X}}s(n,y)\ \ \ \mbox{where }\ \ \  s(n,y)=1-\frac{\pi_n(y)}{\pi(y)}\]

and satisfies $d(n)\leq s(n)$.

The following result was established by Aldous and Diaconis (1987) :
\begin{itemize}
 \item If $T$ is a strong stationary time, for all $n\in \N$, $s(n)\leq \Pr(T>n)$.
\item Conversly, there exists a strong stationary time $T$ such that there is equality for all $n\in \N$.
\end{itemize}
\bigskip

The construction of the stochastically optimal strong stationary time is the following :

$\gamma_n(x)=(s(n-1)-s(n))\pi(x)$, so as before $\ds\sum_{x\in {\mathbb X}}\gamma_n(x)=s(n-1)-s(n) $

With the  method of proposition \ref{prop:defpsi} with $\psi_0(x)=(1-s(0))\frac{\pi(x)}{\pi_0(x)}$ on the support of $\pi_0$, we obtain $\Pr(X_n=x,T=n)=(s(n-1)-s(n))\pi(x)$. So $T$ is an optimal strong stationary time.

\erque

 \brque
We see in the preceding remark that in the cases of separation and total variation, we can construct a randomized stopping time $T$ such that $d(n)=\Pr(T>n)$ for all $n\in \N$, the construction depends of the choice of $\gamma_n$ such that $\ds\sum_{x\in {\mathbb X}}\gamma_n(x)=d(n-1)-d(n) $.

For example, if for the total variation distance, we take $\gamma_n(x)=(d(n-1)-d(n))\pi(x)$ as for separation, we can prove that it does not exist $\psi_n$ with values in $[0,1]$ which satifies $\Delta_n\psi_n=\gamma_n$.

We can deal with others distances and the first condition is to have $n\rightarrow d(n)$ decreasing.

For any convex function $f:\R_+\rightarrow \R$, one may define the $f$-divergence distance on ${\cal P}(E)$ by $d_f(\mu,\pi)=\ds\sum_{x\in {\mathbb X}}\pi(x)f(\frac{\mu(x)}{\pi(x)})$.

These measures are studied in \cite{Liese87} by Liese-Vajda.

We suppose that $f(1)=0$, so $d_f(\mu,\pi)\geq 0$.

Some examples :
\begin{itemize}
 \item The total variation distance given by $f(x)=\demi|x-1|$.
\item The relative entropy given by $f(x)=x\ln (x)$.
\item The Hellinger distance given by $f(x)=\demi (1-\sqrt{x})^2$.
\item The $\chi^2$ distance given by $f(x)=\demi (1-x)^2$.
\end{itemize}
We note $d_f(n)=d_f(\pi_n,\pi)$.

Let $\gamma_n:{\mathbb X}\rightarrow \R$ be defined by $\gamma_n(x)=(\pi f(\frac{\pi_{n-1}}{\pi}))K-f(\frac{\pi_{n}}{\pi})\pi$. So $d_f(n-1)-d_f(n)=\ds\sum_{x\in {\mathbb X}}\gamma_n(x)$. The convexity of $f$ implies that $\gamma_n(x)\geq 0$.

In all generality, $d_f(n)$ can be superior to $1$ and so it will be impossible to find a randomized stopping time $T$ such that for all $n\geq 0$, $\Pr(T>n)=d(n)$.

So a constant $c>0$ and a function $\psi_0:{\mathbb X}\rightarrow [0,1]$ have to be found such that for all $n\geq 0$, $\Pr(X_n=x,T=n)=\frac{\gamma_n(x)}{c}=\psi_n(x)\Delta_n(x)$.

If we want $\Pr(T<+\infty)=1$, it is necessary to choose $\psi_0$ such that
\[1=\Pr(T=0)+\ds\sum_{n\geq 1}\Pr(T=n)=\sum_{x\in {\mathbb X}}\psi_0(x)\pi_0(x)+\frac{1}{c}d_f(0)\]

By mimicking the preceding calculus, for all $n\geq 1$, 
\[\Delta_n(x)=\pi_n(x)-\left[(\psi_0\pi_0)+\frac{1}{c}f(\frac{\pi_0}{\pi})\pi \right]K^n(x)+\frac{1}{c}f(\frac{\pi_n}{\pi})\pi(x)+\frac{1}{c}\gamma_n(x)\]

So the condition $\Delta_n\geq \frac{1}{c}\gamma_n$ implies that 
\[ f(\frac{\pi_n}{\pi})\pi\geq -c\pi_n+\left[c(\psi_0\pi_0)+f(\frac{\pi_0}{\pi})\pi \right]K^n\]

In the case where $\pi_0=\delta_{x_0}$, this condition becomes $f(\frac{\pi_n}{\pi})\geq f(0)(1-\frac{\pi_n}{\pi})$.

But as $f(1)=0$ and $f$ is convex, if $\frac{\pi_n(x)}{\pi(x)}\leq 1$, $f(\frac{\pi_n(x)}{\pi(x)})\leq (1-\frac{\pi_n(x)}{\pi(x)})f(0)$.

This implies that  $f/[0,1]$  is affine.

In our examples, the only one is the total variation distance.

\erque

\subsection{Construction of an optimal stopping time in continuous  time}
Let $Q=(q(x,y))_{(x,y)\in {\mathbb X}\times {\mathbb X}}$ be the generator of an irreducible positive recurrent  Markov chain on the finite space ${\mathbb X}$.

As before, $\pi_0$ is a probability measure on ${\mathbb X}$. $(X_t)_{t\geq 0}$ is a Markov chain with initial distribution $\pi_0$ and generator $Q$.

We write $\pi_t=\pi_0\exp(tQ)$ for the distribution of $X_t$.

Let $\pi$ be the unique stationary distribution. $\pi$ is the unique distribution satisfying $\pi Q=0$. 

We write $d(t)=\|\pi_t-\pi\|_{TV}$.

Let $(J_t)_{t\geq 0}$ be defined by
\begin{equation}\label{eq:defJ}
J_t=\left\lbrace \begin{array}{ll}
                  0&\mbox{if } \ \exists s\in [0,t],\ \pi_s(X_s)\leq \pi(X_s)\\
                  (1-\psi_0(X_0))exp\left( \int_0^t\frac{(\pi_s\wedge \pi)Q(X_s)}{\pi_s(X_s)-\pi(X_s)}ds\right) &\mbox{else}
                 \end{array}
\right.  
\end{equation} 
where $\psi_0$ is defined in (\ref{psi}).
Let $U$ be a uniformly distributed on  $[0,1]$ random variable independent of $(X_t)_{t\geq 0}$.

A randomized stopping time $T$ is defined by $T=\inf\{t\geq 0/\ U\geq J_t\}$.

\bprop\label{prop:weakcontcoup}
\begin{equation}
 \forall t\geq 0, d(t)=\Pr(T>t)
\end{equation} 
\begin{equation}
 \forall t\geq 0, \Pr(T\leq t, X_t=x)=(\pi_t\wedge \pi)(x)
\end{equation} 
\begin{equation}
 \forall t\geq 0, \Pr(T\leq t, X_T=x)=(\pi_t\wedge \pi)(x)-\int_0^t(\pi_s\wedge \pi)Q(x)ds
\end{equation} 
The distribution $\nu$ of $(T,X_T)$ on $(\R_+\cup \{+\infty\})\times {\mathbb X}$ is given by
\begin{equation}
 \nu(ds,x)=(\pi_0\wedge\pi)(x)\delta_0(s)+(\un_{\pi_s(x)<\pi(x)}\pi_sQ(x)-(\pi_s\wedge\pi)Q(x))ds
\end{equation} 
\eprop

The proof of the proposition (\ref{prop:weakcontcoup}) can be found in part ({\bf 5.2}) of the appendix.

\subsection{Some remarks about weak and maximal coupling}

In the preceding sections, the construction of an optimal randomized stopping time $T$ such that for all $n$, $\Pr(T>n)=d(n)$ is given. It is easy to see that by the same way, a randomized stopping time $\widetilde{T}$ for a Markov chain $(Y_n)_{n\geq 0}$ whose initial distribution given by $\pi$, can be defined such that $\Pr(\widetilde{T}>n)=d(n)$ and $\Pr(X_n=x,T=n)=\Pr(Y_n=x,\widetilde{T}=n)=\gamma_n(x)$.

So $((X_n)_{n\geq 0},T)$ and $((Y_n)_{n\geq 0},\widetilde{T})$ is a weak optimal coupling as defined in \cite{Lindvall}. Then Lindvall prove that it exists  $(\widehat{X},\widehat{Y})$ a coupling of $(X,Y)$ such that the coupling time $\widehat{T}$ of $(\widehat{X},\widehat{Y})$ has the same distribution as $T$.

Furthermore, the coupling process $(\widehat{X},\widehat{Y})$ can be chosen such that it is a time inhomogeneous Markov chain. The transitions can be computed but are quite complicated. It is not necessarily a Markovian coupling. By definition a Markovian coupling $((X_n,Y_n)_{n\geq 0}$ of $(\pi_0,K),(\pi,K)$ satisfies the following property :
\[\forall N, \ \mbox{under } \Pr(\ \ |(X_k,Y_k)_{0\leq k\leq N}),\ \ (X_{N+n},Y_{N+n})_{n\geq 0}\mbox{ is a coupling of }\ (\pi_N,K),(\pi,K)\]
This is sometimes called causal or co-adapted coupling. These ones do not always exist.

Furthermore, it can be proved that $\widetilde{T}=\min(S,\widetilde{S})$ where $S=\inf\{n\geq 0/ \pi_n(\widetilde{X}_n)\leq \pi(\widetilde{X}_n)\}$ and $
\widetilde{S}=\inf\{n\geq 0/ \pi_n(\widetilde{Y}_n)\geq \pi(\widetilde{Y}_n)$.

Here is an example of what is obtained for the two-state chain in continuous time.

${\mathbb X}=\{0,1\}$, $Q=\left(\begin{array}{cc}
                          -\lambda&\lambda\\
\mu&-\mu
                         \end{array}
 \right) $, $\theta=\lambda+\mu$ and $\pi_0=p\delta_0+q\delta_1$.

The stationary distribution $\pi$ is given by $\pi=\frac{\mu}{\theta}\delta_0+\frac{\lambda}{\theta}\delta_1$ and suppose that $p>\frac{\mu}{\theta}$.

We have $\pi_t=(\frac{\mu}{\theta}+\frac{p\lambda-\mu q}{\theta}e^{-t\theta})\delta_0+(\frac{\lambda}{\theta}-\frac{p\lambda-\mu q}{\theta}e^{-t\theta})\delta_1$ and $d(t)=\frac{p\lambda-\mu q}{\theta}e^{-t\theta}$.

The initial law of the coupling is given by what is called the  $\gamma$-coupling of $(\pi_0,\pi)$ in \cite{Lindvall}, 
\[\Pr((X_0,Y_0)=(0,0))=\pi(0),\Pr((X_0,Y_0)=(1,1))=\pi(1),\Pr((X_0,Y_0)=(0,1))=p-\pi(0)\]

We find in this case an homogeneous Markov chain on ${\mathbb X}\times {\mathbb X}$ with generator given by
\[q((0,1),(1,1))=\lambda,q((0,1),(0,0))=\mu,q((1,1),(0,0))=\mu,q((0,0),(1,1))=\lambda\]

The distribution of the coupling time $T^*$ is $(q+\frac{\mu}{\theta})\delta_0+(p-\frac{\mu}{\theta}){\cal E}(\theta)$.

In \cite{Fill1991}, Fill finds that in this example, the optimal coupling time for separation is given by an exponential random variable with parameter $\theta$. 

\section{Halting set for total variation distance}
\subsection{Proof of theorem \ref{theo:condsufhalting}, existence of halting set}
By the preceding constructions, for all $t\geq 0$, $d(t)=\Pr(T>t)=\Pr(U<J_t)$. By definition of $J_t$ (\ref{eq:defJ}), (\ref{Jn}), $\{U<J_t\}\subset \{T_{\cal M}>t\}$. So $T\leq T_{\cal M}$.

\brque\label{re:eigenvalue}
How can we find a halting state ?

Suppose that our Markov chain is a reversible aperiodic Markov chain.

Let $|{\mathbb X}|=N+1$ and let $1=\beta_0>\beta_1\geq \cdots\geq \beta_N>-1$ be the eigenvalues of $K$ with $L^2(\pi)$-normalized eigenvectors $V_0,\ldots,V_N$. The spectral decomposition gives:
\[\frac{\pi_n(x)}{\pi(x)}-1=\ds\sum_{k=1}^N\beta_k^n(\sum_{y\in {\mathbb X}}\pi_0(y)V_k(y)V_k(x))\]

Denote by $\rho=\max(|\beta_i|,1\leq i\leq N)$. So if an halting state $x^*$ exists, it must satisfy
\begin{equation} \label{eq:condop}
 \ds\sum_{k/|\beta_k|=\rho}\sum_{y\in {\mathbb X}}\pi_0(y)V_k(y)V_k(x_*)\leq 0\ \mbox{and } \ds\sum_{k/|\beta_k|=\rho}\mbox{sign} (\beta_k)\sum_{y\in {\mathbb X}}\pi_0(y)V_k(y)V_k(x_*)\leq 0
\end{equation} 
\erque

Of course, the set $\cal M$ can be empty. Property of monotonicity can imply the existence of halting state, as proved in the following proposition.

Let $\prec$ be a partial order on ${\mathbb X}$. 

A Markov chain on ${\mathbb X}$ with transition kernel given by $K=(k(x,y))_{x,y\in {\mathbb X}}$ is monotone if for all $x \prec y$, the probability $K(x,\cdot)$ is stochastically smaller than $K(y,\cdot)$. That means that for an increasing function $f$, $Kf$ is increasing.
\bprop\label{prop:monotony}
Let $K=(k(x,y))_{x,y\in {\mathbb X}}$ a monotone kernel for a partial order on ${\mathbb X}$. Suppose furthermore that there exists a smallest element denoted by $\mathbf{0}$ and a largest element denoted by  $\mathbf{1}$.

Then $\mathbf{1}$ is a halting state for the Markov chain started at $\mathbf{0}$.
\eprop

\begin{proof}
 $\frac{\pi_0}{\pi}$ is decreasing. We prove inductively that for all $n\geq 0$, $\frac{\pi_n}{\pi}$ is decreasing.
 
 It comes from the relation $\frac{\pi_{n+1}}{\pi}=K(\frac{\pi_n}{\pi})$.
 
 So if $\frac{\pi_n}{\pi}(\mathbf{1})>1$, for all $x\in {\mathbb X}$, $\pi_n(x)>\pi(x)$ that is false.
\end{proof}
\brque
In the continuous case, the hypothesis is changed to the following, for all $t\geq 0$, the kernel $P_t=e^{tQ}$ is monotone.

\erque

\subsection{One example: The riffle shuffle }
The well-known rifle-shuffle is a method of shuffling cards. Its mathematical description was made by Gilbert and Shannon (see \cite{Gilbert}) and independently by Reeds \cite{Reeds}. A sharp mathematical analysis for the 
riffle shuffling was carried out by Bayer and Diaconis (1992) \cite{Bayerdiaconis}.

Denote by $(\sigma_n)_{n\geq 0}$ this Markov chain on the symmetric group ${\mathfrak S}_N$. $\sigma_0$ is the identity and the stationary measure is given by the uniform one.

Let $f(\sigma)$ be the number of rising sequence of $\sigma$. In \cite{Bayerdiaconis}, corollary 2. it is proved that $(f(\sigma_n))_{n\geq 0}$ is a Markov chain on ${\mathbb X}=\{1,\ldots , N\}$. Its stationary probability is given by $\nu(r)=\frac{A_{N,r}}{N!}$ where $A_{N,r}$, the Eulerian number,  is the number of permutations with $r$ rising sequences.

If we denote by $\nu_n$ the law of $f(\sigma_n)$, we have $\|\pi_n-\pi\|_{TV}=\|\nu_n-\nu\|_{TV}$.

In \cite{Bayerdiaconis}, theorem 3, it is proved that $\Pr(\sigma_n=\sigma)=\frac{\binom{N+2^n-r}{N}}{2^{nN}}$ if $f(\sigma)=r$.

So a halting state for total variation for the chain $(f(\sigma_n))_{n\geq 0}$ is a state $r$ such that for all $n\geq 0$, $\frac{\binom{N+2^n-r}{N}}{2^{nN}}\leq\frac{1}{N!}$.

We can see easily that $[\frac{N+1}{2}]$ is a halting state for total variation.

Furthermore $(f(\sigma_n)-1)_{n\geq 0}$ is a Markov chain on $\{0,\ldots ,N-1\}$ which is symmetric, i.e. $k(x,y)=k(N-1-x,N-1-y)$. We shall see later that for symmetric birth--and-death chain, the middle element is often a halting state.

If ${\mathfrak M}$ is the set of permutations with more than $ [\frac{N+1}{2}]$ rising sequence then we have $\|\pi_n-\pi\|_{TV}\leq \Pr(T_{\mathfrak M}>n)$.

We deduce using Theorem 4 of  (\cite{Bayerdiaconis}) that $\Pr(T_{\mathfrak M}>n)\geq 1-2\Phi(\frac{-1}{4c\sqrt{3}})+O_c(\frac{1}{N^{\frac{1}{4}}})$ where $0<c$, $n=\log_2(N^{\frac{3}{2}}c)$ and $\Phi(x)=\ds\int_{-\infty}^xe^{-\frac{1}{t^2}}\frac{dt}{\sqrt{2\pi}}$.

So for example, $\varliminf\frac{\Es[T_{\mathfrak M}]}{\log_2(N^{\frac{3}{2}})}\geq 1$.

\subsection{Applications to the bounding of the total variation distance and to problems of cutoffs}
To prove the proposition \ref{prop:cutoffgeneral}, the following results found in the book of Aldous and Fill \cite{AldousFill} are used.

They are given in the discrete time.

Let $P$ be an irreducible transition matrix on a finite space ${\mathbb X}$ with a reversible distribution $\pi$.

Let $A\subset {\mathbb X}$ and $P^A=P/A^c\times A^c$. Let $\lambda_A$ be the largest eigenvalue of $P^A$ and $\lambda _1$ the second eigenvalue of $P$. $P^A$ is supposed to be irreducible and aperiodic, then:

  \begin{equation}\label{propertyeigenvalue}1-\lambda_A\geq \pi(A) (1-\lambda_1)
\end{equation}
The proof can be found in \cite{AldousFill}, Theorem 33, Corollary 34, Chapter 3 Reversible Markov Chains.

 For all $n\geq 0$, 
\begin{equation}\label{propertyTAlambdaA}\Pr_{\pi}(T_A>n)\leq \pi(A^c)\lambda_A^n\end{equation}

The proof can be found in \cite{AldousFill}, Proposition 21, Chapter 3 Reversible Markov Chains. This is an application of the Perron-Froebenius theorem to $P^A$.

These results are also used in \cite{DingLubetzkyPeres} and \cite{BasuHermonPeres}.

The third useful result, seen in \cite{BasuHermonPeres} is the following :
\begin{equation}\label{propertyPeres}\forall x\in {\mathbb X}, \forall n,m,\ \Pr_x(T_A>n+m)\leq \|\delta_xP^n-\pi\|_{TV}+\Pr_{\pi}(T_A>m)
\end{equation}
We prove this result for clarity:
\begin{proof}
 Let $(X_n)_{n\geq 0}$ be $(\delta_x,P)$ Markov chain. Let $n\geq 0$.
 
 We can find a random variable $U$ independent of $(X_n)_{n\geq 0}$ and a function of $X_n$ and $U$ denoted by $Y_n$ such that the law of $Y_n$ is $\pi$ and $\|\delta_xP^n-\pi\|_{TV}=\Pr(X_n\not=Y_n)$. It is a standard argument of coupling.
 
 We have  $\Pr_x(T_A>n+m)\leq \Pr(X_n\not=Y_n)+\Pr(T_A>n+m,\ X_n=Y_n)$. Using the Markov property, we have
 
 \[\begin{array}{ll}\Pr(T_A>n+m,\ X_n=Y_n)&=\Es[\un_{T_A>n, X_n=Y_n}\Pr(T_A(X)>n+m|\sigma(X_k,k\leq n)\vee \sigma(U))]\\&=\Es[\un_{T_A>n, X_n=Y_n}\Pr_{X_n}(T_A(X)>m)]\end{array}\]
 
 So $\Pr_x(T_A>n+m)\leq \|\delta_xP^n-\pi\|_{TV}+\Es[\Pr_{Y_n}(T_A(X)>m)]$.
\end{proof}

{\bf Proof of proposition \ref{prop:cutoffgeneral}}

The proof is similar to these found in the paper of Ding, Lubetzky, Peres (\cite{DingLubetzkyPeres}) and of  Diaconis and Saloff-Coste (\cite{DiaconisSaloff2}) or in \cite{BasuHermonPeres}.

\begin{description}
 \item[Bounding the total variation from above]

For all $n\geq 0$, $d^{(N)}(n)\leq \Pr(T_{A^{(N)}}>n)$. So for all $\gamma>0$,
\[d^{(N)}(\lceil\Es[T_{A^{(N)}}]+\gamma\sigma(T_{A^{(N)}})\rceil)\leq \Pr(T_{A^{(N)}}-\Es[T_{A^{(N)}}]>\gamma \sigma(T_{A^{(N)}}))\leq\frac{1}{1+\gamma^2}\]

\item[Bounding the total variation from below]
By (\ref{propertyPeres}),
\[d^{(N)}(n)\geq \Pr(T_{A^{(N)}}>n+m)-\Pr_{\pi}(T_{A^{(N)}}>m)\]

Let $0<\epsilon<1$ and $n=[(1-\epsilon)\Es[T_{A^{(N)}}]]$ and $m=[\frac{\epsilon}{2}\Es[T_{A^{(N)}}]]$. Then by (\ref{propertyTAlambdaA}) , 
\[\begin{array}{ll}
   \Pr_{\pi}(T_{A^{(N)}}>m)&\leq (1-\pi(A^{(N)}))e^{-m(1-\lambda_{A^{(N)}})}\\
   &\leq ee^{-\frac{\epsilon}{2}(1-\lambda_{A^{(N)}})\Es[T_{A^{(N)}}]}
  \end{array}\]
So 
\[\begin{array}{ll}
 \Pr(T_{A^{(N)}}>n+m)&\geq  \Pr(T_{A^{(N)}}- \Es[T_{A^{(N)}}]\geq -\frac{\epsilon}{2}\Es[T_{A^{(N)}}])\\
 &\geq 1-\frac{\sigma(T_{A^{(N)}})^2}{(\frac{\epsilon}{2}\Es[T_{A^{(N)}}])^2}
  \end{array}\]
  
  D'o\`u $d^{(N)}([(1-\epsilon)\Es[T_{A^{(N)}}]])\geq 1-\big(\frac{\sigma(T_{A^{(N)}})^2}{\frac{\epsilon}{2}\Es[T_{A^{(N)}}]}\big)^2-ee^{-\frac{\epsilon}{2}(1-\lambda_{A^{(N)}})\Es[T_{A^{(N)}}]}$
  
  If $n=[\Es[T_{A^{(N)}}]-\gamma\sigma(T_{A^{(N)}})]$, we obtain using the same inequalities,
  \[d^{(N)}(n)\geq 1-\frac{4}{\gamma^2}-e^{-\frac{\gamma}{2}(1-\lambda_{A^{(N)}})\sigma(T_{A^{(N)}})}\]

\end{description}

\brque
If $\varliminf\pi(A^{(N)})>0$, the condition $(1-\lambda_{A^{(N)}})\Es[T_{A^{(N)}}]\rightarrow +\infty$ can be changed to $(1-\lambda_1^{(N)})\Es[T_{A^{(N)}}]\rightarrow +\infty$ by the property (\ref{propertyeigenvalue}).
\erque

\section{Birth-and-death chain on ${\mathbb X}^{(N)}=\{0,\ldots,N\} $ started at $0$}

\subsection{General results}
For every  birth-and-death chain on $\{0,\cdots,N\}$, we will denote  
\[\begin{array}{ll}
 p_x=k(x,x+1)&\mbox{ if } 0\leq x\leq N-1\\
q_x=k(x,x-1)&\mbox{ if } 1\leq x\leq N\\ 
r_x=k(x,x)& \mbox{ if } 0\leq x\leq N
  \end{array}\]
We have for all $0\leq x\leq N$, $p_x+r_x+q_x=1$. We suppose the Markov chain irreducible, so for all $0\leq x\leq N-1$, $p_x>0$ and for all $1\leq x\leq N$, $q_x>0$.

The chain is started at $0$. For $x,y\in {\mathbb X}^{(N)}$, $T_x$ is the hitting time of the state $x$ and $\tau_{x,y}$ is the hitting time of $y$ by the chain started at $x$.

To use proposition (\ref{prop:cutoffgeneral}), it is necessary to know the spectrum of the transition matrix. So below, is given a result without eigenvalues.
\bprop\label{prop:cutoffbirth}
Suppose that  $x^*_N$ is a halting state for total variation for the chain started at $0$ such that $\Es[T_{x^*_N}]$ goes to $+\infty$ when $N$ goes to $+\infty$ and such that $\sigma(T_{x^*_N})=o(\Es[T_{x^*_N}])$.

Suppose also that there exists $y<x^*_N$ such that $\Es[T_{x^*_N}]\sim \Es[T_{y}]$ and $\pi^{(N)}(\{0,\ldots,y\})\rightarrow 0$;

Then the chain has a  cutoff at time $\Es[T_{x^*_N}]$ for total variation.

If furthermore $\Es[\tau_{y,x^*_N}]=O(\sigma(T_{x^*_N}))$, then there is a   cutoff at time $\Es[T_{x^*_N}]$ with window $\sigma(T_{x^*_N})$ for total variation.

\eprop

\begin{proof}
 To bound the total variation distance from below, the following inequality is used :
 
For all $n\geq 0$, $d^{(N)}(n)\geq \Pr(T_y\geq n)-\pi^{(N)}(\{0,\ldots,y\})$.
 
If $n=[(1-\epsilon)\Es[T_{x^*_N}]]$,
\[\begin{array}{ll}\Pr(T_y\geq n)&=1-\Pr\left( \Es[T_y]-T_y>\Es[T_{x^*_N}](\epsilon-(1-\frac{\Es[T_y]}{\Es[T_{x^*_N}}))\right)\\
 &\geq 1-\frac{\sigma(T_{y})^2}{\Es[T_{x^*_N}]^2(\epsilon-(1-\frac{\Es[T_y]}{\Es[T_{x^*_N}]}))^2}\\
&\geq 1-\frac{\sigma(T_{x^*_N})^2}{\Es[T_{x^*_N}]^2(\epsilon-(1-\frac{\Es[T_y]}{\Es[T_{x^*_N}]}))^2}
\end{array}\]

So $d^{(N)}(n)\geq 1-\frac{\sigma(T_{x^*_N})^2}{\Es[T_{x^*_N}]^2(\epsilon-(1-\frac{\Es[T_y]}{\Es[T_{x^*_N}]}))^2}-\pi^{(N)}(\{0,\ldots,y\})$.

According to the hypothesis, it is clear that $\ds\lim_{N\rightarrow +\infty}d^{(N)}([(1-\epsilon)\Es[T_{x^*_N}]])=1$.

If furthermore $\Es[\tau_{y,x^*_N}]=O(\sigma(T_{x^*_N}))$, if $\gamma$ big enough,
\begin{equation}\label{eq:above}
 d^{(N)}(\Es[T_{x^*_N}]-\gamma\sigma(T_{x^*_N}))\geq 1-\frac{1}{(\gamma-\frac{\Es[\tau_{y,x^*_N}]}{\sigma(T_{x^*_N})})^2}-\pi^{(N)}(\{0,\ldots,y\}
\end{equation}

\end{proof}

There are several methods to calculate $\Es[T_x]$ and $\sigma(T_x)$ for a birth-and-death chain. The first one uses the spectrum with the following formulae, see \cite{KarlinMcGregor} :
\begin{equation}\label{hittingbyeigen}
\Es[T_x]=\ds\sum_{k=0}^{x-1}\frac{1}{1-\alpha_k}\ ,\ \sigma(T_x)^2=\ds\sum_{k=0}^{x-1}\frac{\alpha_k}{(1-\alpha_k)^2}
 \end{equation}
where $\alpha_0>\ldots>\alpha_{x-1}$ are the eigenvalues of $P_{/\{0,\ldots,x-1\}^2}$. This implies the following inequality:
\begin{equation}\label{propertysigmaesp}
 \sigma(T_x)^2\leq \frac{1}{1-\alpha_0}\Es[T_x]
\end{equation}

The second one uses the following formulae :
\begin{equation}
 \Es[T_x]=\ds\sum_{y=0}^{x-1}\Es[\tau_{y,y+1}]\ \mbox{ and }\ \sigma(T_x)^2=\ds\sum_{y=0}^{x-1}\sigma(\tau_{y,y+1})^2
\end{equation}
\begin{equation}\label{eq:hittingtimes}
\Es[\tau_{x,x+1}]=\frac{\pi(\{0,\ldots,x\})}{p_x\pi(x)}\ ,\ \sigma(\tau_{x,x+1})^2=\frac{2}{\pi(x)p_x}\ds\sum_{y=0}^{x-1}\frac{\pi(\{0,\ldots,y\})^2}{p_y\pi(y )}+\Es[\tau_{x,x+1}]^2-\Es[\tau_{x,x+1}]
\end{equation}
Sometimes, it is  possible to calculate $\Es[T_x]$ by using martingale.

\subsection{Monotone birth-and-death chain on ${\mathbb X}^{(N)}=\{0,\ldots,N\} $ started at $0$}

A birth-and-death chain is monotone if for all $0\leq x\leq N-1$, $p_x+q_{x+1}\leq 1$.

In this case if $\pi_0=\delta_0$ and as usual, $\pi_n=\pi_0K^n$, then for all $n\geq 0$, $x\mapsto \frac{\pi_n(x)}{\pi(x)}$ is decreasing and by consequence $N$ is a halting state for total variation as proved in proposition (\ref{prop:monotony}).

An example of monotone birth-and-death chain is given by birth-and-death chain with positive spectrum and in particular by lazy chain which are chain with $r_x\geq\demi$ for all $x$.

\brque
With the notation of the section \ref{sec:optstopti}, it can be proved that in the case of a monotone birth-and-death chain on $\{0,\cdots,N\} $ started at $0$, the support of the distribution of $X_T$ is  $\{0,\ldots,x*\}$ where $x*$ is the smallest halting state for total variation.

Indeed, $\Pr(X_T=x)=\displaystyle\sum_{n=0}^{+\infty}\gamma_n(x)$.

So $\Pr(X_T=0)\geq \gamma_0(0)>0$, $\Pr(X_T=1)\geq \gamma_1(1)=p_0\wedge \pi(1)-p_0\pi(0)>0$

By iteration and by using monotonicity of $x\mapsto \frac{\pi_n(x)}{\pi(x)}$,  
\[\forall n\geq 0,\ \gamma_n(x)=0\Leftrightarrow \forall n\geq 0,\ \pi_{n}(x-1)\leq \pi(x-1)\]
By consequence, the support of the distribution of $X_T$ is $\{0,\ldots,x^*\}$ where for all $n\geq 0$, $\pi_n(x^*)\leq\pi(x^*)$ and there exists $n\geq 0$ with $\pi_n(x^*-1)>\pi(x^*-1)$.

So $x^*$ is the smallest halting state for total variation.

The proof is the same in continuous time using (\ref{prop:weakcontcoup}).
\erque

There is an interesting property satisfied by monotone birth-and-death chain. With the  notations of (\ref{re:eigenvalue})

\blemm \label{le:gap=second}

If $P$ is a monotone irreducible birth-and-death matrix with eigenvalues given by  $1=\beta_0>\beta_1\geq \cdots\geq \beta_N>-1$, then $|\beta_N|\leq \beta_1$.
\elemm 
\begin{proof}

The spectrum of an irreducible birth-and-death chain satisfies $1=\beta_0>\beta_1>\cdots>\beta_N$, so $\rho\in\{\beta_1,|\beta_N|\}$

We have seen by monotony that for all $n\geq 0$, $\frac{\pi_n(N)}{\pi(N)}\leq 1$. But
\[\begin{array}{ll}\frac{\pi_n(N)}{\pi(N)}-1&=\ds\sum_{k=1}^N\beta_k^nV_k(0)V_k(N)\\
 &=\beta_1^n  V_1(0)V_1(N)+\beta_N^nV_N(0)V_N(N)+\ds\sum_{k=2}^{N-1}\beta_k^nV_k(0)V_k(N)\\
  \end{array}
\]

But by using property (\ref{eq:signeigenvector}) of the appendix on the eigenvectors of irreducible birth-and-death chain, we have that $V_1(0)V_1(N)<0$ and $(-1)^NV_N(0)V_N(N)>0$.

So if $|\beta_N|>\beta_1$, $\frac{\pi_n(N)}{\pi(N)}-1\sim \beta_N^nV_N(0)V_N(N)$ and we have a contradiction with the fact that for all $n\geq 0$, $\frac{\pi_n(N)}{\pi(N)}-1\leq 0$.

By consequence $|\beta_N|\leq \beta_1$.
\end{proof}
\bexem {\bf Metropolis chains}

Let $\pi$ be a probability on ${\mathbb X}^{(N)}=\{0,\ldots,N\} $ with $\pi(x)>0$ for all $x\in {\mathbb X}^{(N)}$. Use the Metropolis algorithm with base chain the simple symmetric random walk to obtain a birth-and- death chain with stationary measure $\pi$ (see e.g.(\cite{DiaconisSaloff})). By construction, it is a monotone chain and so $N$ is a halting state. 
\eexem

\bexem \label{exem:srw}{\bf Simple random walk}

Consider the birth-and-death chain on ${\mathbb X}^{(N)}=\{0,\ldots,N\} $ started at $0$ with  $p_x=p$, $q_x=q$, $r_0=q$ and $r_N=p$. We suppose $0\leq q<p\leq 1$ and $p+q=1$. It is a monotone chain, so $N$ is a halting state.

The spectrum of this kernel is known, see \cite{Feller} for example.

The eigenvalues are $\beta_{k}=2\sqrt{pq}\cos(\frac{\pi k}{N+1})$ for $1\leq k\leq N$ and $V_k$ is proportional to $\left( \begin{array}{c}x_{k,0}\\ \vdots\\ x_{k,N}\end{array}\right) $ where $x_{k,l}=\sqrt{p}\sin(\frac{\pi k(l+1)}{N+1})-\sqrt{q}\sin(\frac{\pi kl)}{N+1})$.

Some calculus prove that (\ref{eq:condop}) is satisfied if $\cot (\frac{\pi x_*}{N+1})\leq -(1-2q)\cot(\frac{\pi}{N+1})$.

If $x_N$ is the smallest integer such that  $x_N>N+1-\frac{N+1}{\pi}\cot^{-1}((1-2q)\cot(\frac{\pi}{N+1}))$, we have $x_N\simeq N+1-\frac{1}{1-2q}$ which is close to $N$ for $N$ large enough.

It can be proved, using methods of martingales that for all $0\leq x\leq N$,
\[\Es[T_x]=\frac{x}{p-q}-\frac{q}{(p-q)^2}(1-(\frac{q}{p})^x)\ \ \ \sigma^2(T_x)=Ax+B+C(\frac{q}{p})^x+D(\frac{q}{p})^{2x}+Ex(\frac{q}{p})^x\]
with $A=\frac{4pq}{(p-q)^3}$ and $B,C,D,E$ some constants depending of $p$.

So proposition(\ref{prop:cutoffbirth}) implies that there is a cutoff at time $t_N$ with a window $\sigma_N$, where $t_N=\frac{N}{p-q}$ and $\sigma_N=\sqrt{ N}$.

To prove this, one can take $y=N-\sqrt{N}$.

If $p=q=\demi$, we have $\Es[T_x]=x(x+1)$ and $\sigma^2(T_x)=\frac{x(x+1)(2x^2+2x-1)}{3}$.

It is well known that in this case, there is no cutoff.

We can prove that in this case $x^*_N=[\frac{N}{2}]+1$ is the smallest halting state for total variation for the chain started at $0$, but the condition  $\sigma(T_{x^*_N})=o(\Es[T_{x^*_N}])$ is not satisfied.

\brque
We can suppose that $p$ depends of $N$. 

As we have seen above, the smallest  halting state for total variation for the chain started at $0$ is bigger than $N+1-\frac{N+1}{\pi}\cot^{-1}((1-2q)\frac{\cos(\frac{\pi}{N+1})}{\sin(\frac{\pi}{N+1})})$.

Let $\delta>0$ and $p_N=\demi+\frac{\delta}{N}$ and $q_N=\demi-\frac{\delta}{N}$.

So the smallest  halting state for total variation for the chain started at $0$ belongs to $\{[\frac{N}{2}],\ldots, N\}$. Denote it by $x_N^*$. As in the case above, the condition $\sigma(T_{x^*_N})=o(\Es[T_{x^*_N}])$ is not satisfied.

\erque

\brque
If we suppose that $p_N=\demi+\epsilon_N$ and $q_N=\demi-\epsilon_N$ with $\epsilon_N>0, \epsilon_N=o(1)$ and $N\epsilon_N\rightarrow +\infty$.

We have $\sigma^2(T_N)\sim \frac{N}{2\epsilon_N^3}$ and $\Es[T_N]\sim \frac{N}{2\epsilon_N}$. So $\sigma(T_{N})=o(\Es[T_{N}])$.

Furthermore it is  easy to see that $y=N-\sqrt{\frac{N}{4\epsilon_N}}$ satisfies the hypothesis of proposition(\ref{prop:cutoffbirth}), so we have a $\frac{N}{2\epsilon_N}$ cut-off.
\erque
\eexem

\subsection{The case of symmetric birth-and-death chain}
\bdefi A Markov chain on ${\mathbb X}=\{0,\ldots,N\}$ with transitions given by $K=(k(x,y))_{x,y\in {\mathbb X}}$ is symmetric if for all $x\in {\mathbb X}$, $y\in {\mathbb X}$ then $k(N-x,N-y)=k(x,y)$.
\edefi
Two different cases can occur, $N$ odd or $N$ even. The spectral analysis differs in these two cases.

\subsubsection{The case of symmetric birth-and-death chain on $\{0,\ldots,2N+1\}$}

Denote by $P$ the transition matrix of the birth-and-death chain, by $Q$ the restriction of $P$ to $\{0,\ldots,N\}\times\{0,\ldots,N\}$.

Let $L$ be  the $(N+1)\times(N+1)$ matrix define by $L_{i,j}=1$ if $i+j=N$ and  $L_{i,j}=0$ if not.

So $P$ is the tridiagonal matrix given by
$P= \left(\begin{array}{c|c}
  Q

 & \begin{array}{ccc}
       &&\\
       &&\\
       p_N&&\end{array}
\\
 
 \hline
 
 \begin{matrix}
  && q_{N+1}\\
 && \\
&&
\end{matrix}
 &
 
 LQL
 
 \end{array}\right)$

Denote by $Q_1$ the matrix equal to $Q$ unless the entry $Q_1(N,N)$ that is $Q_1(N,N)=r_N+p_N$. So $Q_1$ is a stochastic matrix.

\bprop\label{prop:monotoneodd}
If $P$ is monotone, then $N+1$ is a halting state for total variation distance for the chain started at $0$.
\eprop

\begin{proof}
 If $(X_n)_{n\geq 0}$ is a birth-and-death symmetric chain on $\{0,\ldots,2N+1\}$ with transition matrix given by $P$ then the process $(Z_n)_{n\geq 0}$ given by $Z_n=N+\demi-|X_n-(N+\demi)|$ is a birth-and-death chain on $\{0,\ldots,N\}$ with transition matrix given by $Q_1$. Its stationary probability is given by $\widetilde{\pi}$ with $\widetilde{\pi}(x)=2\pi(x)$.

And if $P$ is monotone then $Q_1$ is monotone too. So for all $n\geq 0$, $\frac{Q_1^n(0,N)}{\widetilde{\pi}(N)}\leq 1$.

Also for all $n\geq 0$, $x\longmapsto\frac{P^n(0,x)}{\pi(x)}$ is decreasing. So
\[\begin{array}{ll}
   2\frac{P^n(0,N+1)}{\pi(N+1)}&\leq\frac{P^n(0,N)}{\pi(N)}+\frac{P^n(0,N+1)}{\pi(N+1)}\\
   &\leq \frac{2}{\widetilde{\pi}(N)}\left(\Pr(X_n=N)+\Pr(X_n=N+1) \right) \\
&\leq  \frac{2}{\widetilde{\pi}(N)}\Pr(Z_n=N)=2\frac{Q_1^n(0,N)}{\widetilde{\pi}(N)}\leq 2
  \end{array}\]

\end{proof}

\bprop
If for all $x\in\{0,\ldots,2N+1\}$, $r_x=0$ then for all $n\geq 0$, for all $k\geq 0$, 
\[\frac{P^n(0,N+1)}{2\pi(N+1)}\leq 1\ \mbox{ and }\frac{P^n(0,N+2)}{2\pi(N+2)}\leq 1\]
\eprop

\begin{proof}
  Denote by $Q_2$ the matrix equal to $Q$ unless the entry $Q_2(N,N)$ that is $Q_2(N,N)=r_N-p_N$. 

We have the following easy lemma :
\blemm\label{lemm:eigenr=0}
If $\lambda$ is an eigenvalue of $Q_1$ associated with the eigenvector $v$, then $\lambda$ is an eigenvalue of $P$ associated with the eigenvector $\begin{pmatrix}
  v\\
Lv                                                                                                                                                                                                                                                                                               \end{pmatrix}$.

If $\lambda$ is an eigenvalue of $Q_2$ associated with the eigenvector $v$, then $\lambda$ is an eigenvalue of $P$ associated with the eigenvector $\begin{pmatrix}
  v\\
-Lv                                                                                                                                                                                                                                                                                               \end{pmatrix}$.

\elemm

Denote by $1=\beta_0>\beta_1>\cdots>\beta_N$ the eigenvalues of $Q_1$.

Suppose $N=2a$, the proof is similar if $N$ is odd.

The formula (\ref{eq:r=0N+1}) of the appendix tells us :
\[P^{2n+1}(0,N+1)=2\pi(N+1)(1+\ds\sum_{k=1}^N\beta_k^{2n+1}(-1)^k\nu_k)\mbox{ with }\nu_k>0 \]
 
If we denote for $0\leq k\leq a$, $\lambda_{2k}=\beta_k$ and for $0\leq k\leq a-1$, $\lambda_{2k+1}=|\beta_{N-k}|$, the sequence $(\lambda_k)_{0\leq k\leq N}$ by property (\ref{eq:r=0eigenvalue}) of the appendix, is decreasing and 
\[P^{2n+1}(0,N+1)=2\pi(N+1)(1+\ds\sum_{k=1}^N\lambda_k^{2n}B_k)\mbox{ where }B_{2k}=(-1)^k\lambda_{2k}\nu_k\mbox{ and }B_{2k+1}=(-1)^{k+1}\lambda_{2k+1}\nu_{N-k}\]
So $B_{2k}B_{2k+1}<0$ and $B_{2k+1}B_{2k+2}>0$. The following lemma gives the result.

\blemm\label{le:le1}
Let $\Gamma_n=\ds\sum_{k=0}^{2d}\alpha^n_kB_k$ where
\begin{equation}\label{le:conGamma}
B_0>0,\ \forall i \ B_{2i}B_{2i+1}<0,\ B_{2i+1}B_{2i+2}>0,\ \alpha_0= 1>\alpha_1>\alpha_2>\ldots>\alpha_{2d}>0
\end{equation}
We suppose that $\Gamma_0=\Gamma_1=\cdots=\Gamma_{d-1}=0$. Then $n\mapsto \Gamma_n$ is stricly increasing on $\{d,d+1,\cdots\}$.

\elemm

The proof is in part ({\bf 5.3}) of the appendix

\end{proof}

\subsubsection{The case of symmetric birth-and-death chain on $\{0,\ldots,2N\}$}

Denote by $Q_+$ the stochastic matrix on $\{0,\ldots,N\}\times\{0,\ldots,N\}$ defined by :
\[\left\lbrace \begin{array}{l}
    Q_+(x,y)=P(x,y)\mbox{ if } 0\leq x,y\leq N-1\\
Q_+(N,N-1)=2p_N,\ Q_+(N,N)=r_N
\end{array}\right.\]

\bprop\label{prop:evenmonotone+}
If $P$ is monotone and $p_{N-1}\leq r_N$ then $N$ is a halting state for total variation distance for the chain started at $0$.
\eprop

\begin{proof}
 If $(X_n)_{n\geq 0}$ is a birth-and-death symmetric chain on $\{0,\ldots,2N\}$ with transition matrix given by $P$ then the process $(Z_n)_{n\geq 0}$ given by $Z_n=N-|X_n-N|$ is a birth-and-death chain on $\{0,\ldots,N\}$ with transition matrix given by $Q_+$. Its stationary probability is given by $\widetilde{\pi}$ with $\widetilde{\pi}(x)=2\pi(x)$ if $x\leq N-1$ and $\widetilde{\pi}(N)=\pi(N)$.

And if $P$ is monotone and $p_{N-1}\leq r_N$ then $Q_+$ is monotone. So for all $n\geq 0$, $\frac{Q_+^n(0,N)}{\widetilde{\pi}(N)}\leq 1$.

But $Q_+^n(0,N)=\Pr(Z_n=N)=\Pr(X_n=N)$ which gives the result.

\end{proof}

\blemm
It exists $\gamma_1,\ldots ,\gamma_N$ such that for all $n\geq 0$, $\frac{P^n(0,N)}{\pi(N)}=1+\ds\sum_{k=1}^N\beta_k^n\gamma_k$ where $1>\beta_1>\cdots >\beta_N$ are the eigenvalues of $Q_+$.
\elemm
\begin{proof}
 If $\lambda$ is an eigenvalue of $Q$ associated with the eigenvector $v$, then $\lambda$ is an eigenvalue of $P$ associated with the eigenvector $\begin{pmatrix}
  v\\
0\\
-Lv                                                                                                                                                                                                                                                                                               \end{pmatrix}$.

If $\lambda$ is an eigenvalue of $Q_+$ associated with the eigenvector $\begin{pmatrix}
  v\\
x                                                                                                                                                                                                                                                                                              \end{pmatrix}$, then $\lambda$ is an eigenvalue of $P$ associated with the eigenvector $\begin{pmatrix}
  v\\
x\\
Lv                                                                                                                                                                                                                                                                                               \end{pmatrix}$.

Spectral properties give the result.

\end{proof}

\bprop\label{prop:positiveeven}
 If all the eigenvalues of $P$ are positive then  $N$ is a halting state for total variation distance for the chain started at $0$.
\eprop

\begin{proof}
 By using the lemma above, for all $n\geq 0$, $\frac{P^n(0,N)}{\pi(N)}=1+\ds\sum_{k=1}^N\beta_k^n\gamma_k$.

So the result comes from the  following lemma (\ref{le:inc}) proved in the appendix with $\Gamma_n=\frac{P^n(0,N)}{\pi(N)}$.

\blemm \label{le:inc}
Let $1=\lambda_0>\lambda_1>\ldots>\lambda_d>0$ and $\Gamma_n=\ds\sum_{i=0}^dA_i\lambda_i^n,A_0\geq 0$

If for $0\leq i\leq d-1$, $\Gamma_i=0$, then $n\rightarrow \Gamma_n$ is increasing and so for all $n\geq 0$, $\Gamma_n\leq A_0$.
\elemm

\end{proof}

\bprop\label{prop:r=0even}
If for all $x\in\{0,\ldots,2N\}$, $r_x=0$ then for all $n\geq 0$,  
\[\mbox{If }2k\geq N,\ \frac{P^{2n}(0,2k)}{2\pi(2k)}\leq 1\]
\[\mbox{ If  }2k+1\geq N,\ \frac{P^{2n+1}(0,2k+1)}{2\pi(2k+1)}\leq 1\]
\eprop

\begin{proof}
 We define two stochastic matrices $Q_e$ and $Q_o$ respectively on  $\{0,\dots,N\}\times\{0,\dots,N\}$ and on $\{0,\dots,N-1\}\times\{0,\dots,N-1\}$ by 
\[Q_e(x,y)=P^2(2x,2y) \ \mbox{ and }\ Q_o(x,y)=P^2(2x+1,2y+1)\]

These stochastic matrices have positive eigenvalues, so they are monotone. Their stationary measures are given respectively by $\pi_e(x)=2\pi(2x)$ and $\pi_o(x)=2\pi(2x+1)$. They are symmetric so we can apply proposition (\ref{prop:monotoneodd}) and proposition(\ref{prop:positiveeven}).

The relations $P^{2n}(0,2x)=Q_e^n(0,x)$ and $P^{2n+1}(0,2x+1)=Q_o^n(0,x)$ give the results.
\end{proof}

\bprop\label{prop:evenrconstant}
 Let $P$ be a transition matrix of a birth-and-death Markov chain on ${\mathbb X}=\{0,\ldots,2N\}$ be given by $P=rI+(1-r)\widetilde{P}$ where for all $x\in {\mathbb X}$, $\widetilde{P}(x,x)=0$. We suppose $0<r<\demi$ and $P$ monotone.
Then for all $n\geq 0$, $P^n(0,N+1)\leq \pi(N+1)$.
\eprop

\begin{proof}
It is an easy consequence of the above proposition (\ref{prop:r=0even}).

For all $n\geq 0$, $P^n=\ds\sum_{k=0}^n\binom{n}{k}r^{n-k}(1-r)^k\widetilde{P}^k$. Let $2x+1\in {\mathbb X}, 2x\geq N$. 

Then $P^{2n}(0,2x+1)=\ds\sum_{k=0}^{n-1}\binom{2n}{2k+1}r^{2n-2k-1}(1-r)^{2k+1}\widetilde{P}^{2k+1}(0,2x+1)$.

So $P^{2n}(0,2x+1)\leq \pi(2x+1)(1-(2r-1)^{2n})\leq \pi(2x+1)$

The same method gives that for $2x\in {\mathbb X}, 2x\geq N$, $P^{2n+1}(0,2x)\leq \pi(2x)$.

Using monotonocity, we have the result.
\end{proof}

\subsubsection{Continuous time symmetric birth-and-death chain on ${\mathbb X}=\{0,\ldots,N\}$}

\bprop
$[\frac{N+1}{2}]$ is a halting state for total variation for the continuous time Markov chain started at $0$.
\eprop

\begin{proof}
 Let $Q$ be the generator of the markov chain. $Q$ is supposed symmetric, $Q(N-x,N-y)=Q(x,y)$.

Let $q\geq 2\max\{q(x),x\in {\mathbb X}\}$ and $P=I+\frac{1}{q}Q$.

$P$ is the transition of a symmetric birth-and-death chain on ${\mathbb X}$ which satisfies $r_x=\Pr(x,x)=1-\frac{q(x)}{q}\geq \demi$. So this chain is lazy and we can apply proposition(\ref{prop:monotoneodd}) and proposition(\ref{prop:positiveeven}) which say that for all $n\geq 0$, $P^n(0,[\frac{N+1}{2}])\leq \pi([\frac{N+1}{2}])$.

But for all $t\geq 0$, $\pi_t([\frac{N+1}{2}])=e^{tQ}(0,[\frac{N+1}{2}])=e^{-qt}\ds\sum_{n=0}^{+\infty}\frac{(tq)^n}{n!}P^n(0,[\frac{N+1}{2}])\leq \pi([\frac{N+1}{2}])$.
\end{proof}

\subsubsection{Applications to cutoffs for symmetric birth-and-death chain}

\bprop\label{prop:cutoff3}
Under the hypothesis that the birth-and-death is  symmetric on $\{0,\ldots,N\}$ such that $x_N^*=[\frac{N}{2}]+1$ is a halting state and  
$(1-\lambda_1)\displaystyle\sum_{k=1}^N\frac{1}{1-\lambda_k}\rightarrow+\infty$, then the chain has a $\Es[T_{x^*_N}]$ cutoff for total variation. 

Furthermore $\Es[T_{x^*_N}]\sim \demi\displaystyle\sum_{k=1}^N\frac{1}{1-\lambda_k}$

If $\liminf (1-\lambda_1)^2\displaystyle\sum_{k=1}^N\frac{\lambda_k}{(1-\lambda_k)^2}>0$, then the chain exhibits a cutoff at time $\Es[T_{x^*_N}]$ with a window $\sigma(T_{x_N^*})$.

Furthermore in this case, the chain exhibits a cutoff at time $\demi\displaystyle\sum_{k=1}^N\frac{1}{1-\lambda_k}$ with a window $\demi\sqrt{\displaystyle\sum_{k=1}^N\frac{\lambda_k}{(1-\lambda_k)^2}}$.

Here $1=\lambda_0>\lambda_1>\cdots >\lambda_N$ are the eigenvalues of the chain.
\eprop

\begin{proof}
 One has to prove that the hypothesis of this proposition imply these of the proposition (\ref{prop:cutoffgeneral}).

% In the case of a symmetric monotone birth-and-death chain on $\{0,\ldots,2N+1\}$, by proposition (\ref{prop:monotoneodd}), $x_N^*=N+1$ is a halting state.

By (\ref{hittingbyeigen}) and using the notation of the appendix, 
\[\Es[T_{N+1}]=\displaystyle\sum_{k=0}^N\frac{1}{1-\gamma_k}\]

But as $\pi(\{N+1,\ldots,2N+1\})=\demi$, by (\ref{propertyeigenvalue}) $\frac{1-\lambda_1}{2}\leq 1-\gamma_0\leq 1-\lambda_1$.

By lemma (\ref{lemm:eigenr=0}) and (\ref{eq:eigenvaluesymm}), for all $0\leq k\leq N$, $\lambda_{2k}=\beta_k$ and $\lambda_{2k+1}=\alpha_k$.

So by (\ref{eq:eigenvaluesymm}), $(1-\lambda_1)\Es[T_{N+1}]$ tends to $+\infty$.

The proof is almost the same in the case of a symmetric birth-and-death chain on $\{0,\ldots,2N\}$. 
\end{proof}

\subsubsection{Two examples}
\paragraph*{The Ehrenfest process }

The Markov kernel of the Ehrenfest chain on $\{0,\ldots,N\}$ is defined by
\[p_x=\frac{N-x}{N+1},\ r_x=\frac{1}{N+1},\ q_x=\frac{x}{N+1}\]

The stationary distribution is given by $\pi(x)=\binom{N}{x}2^{-N}$ and the eigenvalues are $\lambda_k=1-\frac{2k}{N+1},\ 0\leq k\leq N$ with  $L^2(\pi)$-normalized eigenvectors 
\[V_k(x)=\binom{N}{k}^{-\demi}\displaystyle\sum_{i=0}^k(-1)^i\binom{x}{i}\binom{N-x}{k-i}\]

In this case $|\lambda_N|=\lambda_1$. It was proved in lemma (\ref{le:gap=second}) that for a monotone chain, $|\lambda_N|\leq\lambda_1$.

We can search  candidates to be halting states by (\ref{eq:condop}) for the chain started at $0$ :
\[V_1(0)V_1(x)\pm V_N(0)V_N(x)\leq 0\Longleftrightarrow 2x\geq N+1\]

So the smallest candidate is $[\frac{N}{2}]+1$.

If $N$ is odd, we can apply proposition(\ref{prop:monotoneodd}) ans so $\frac{N+1}{2}$ is a halting state for total variation for the chain started at 0.

If $N=2N'$ is even, we haven't $q_{N'+1}\leq r_{N'}$ so we can't apply proposition(\ref{prop:evenmonotone+}).

But $P=rI+(1-r)\widetilde{P}$ with $r=\frac{1}{N+1}$ and $\widetilde{P}(x,x)=0$ for all $x$. So we can apply proposition(\ref{prop:evenrconstant}) and so $N'+1$ is a halting state for total variation for the chain started at 0.

So $x_N=[\frac{N}{2}]+1$ is the smallest halting state for total variation for the chain started at 0.

As $(1-\lambda_1)\displaystyle\sum_{k=1}^N\frac{1}{1-\lambda_k}=\displaystyle\sum_{k=1}^N\frac{1}{k}\rightarrow +\infty$ and $(1-\lambda_1)^2\displaystyle\sum_{k=1}^N\frac{\lambda_k}{(1-\lambda_k)^2}=-\frac{2}{N+1}\displaystyle\sum_{k=1}^N\frac{1}{k}+\displaystyle\sum_{k=1}^N\frac{1}{k^2}\rightarrow \frac{\pi^2}{6}$, proposition (\ref{prop:cutoff3}) gives the existence of a cutoff at time $\demi\displaystyle\sum_{k=1}^N\frac{1}{1-\lambda_k}\sim \frac{N\ln N}{4}$ with a window $N$.

\paragraph*{\bf Bernoulli-Laplace models}\label{ex:bernou}

Consider two urns, the left containing $r$ red balls, the right $N-r$ black balls with $0<2r\leq N$. At each step, a ball is picked uniformly at random in each urn and the two balls are switched. The process is completely determined by the number of red balls in the right urn and this is a birth-and-death chain on ${\mathbb X}=\{0,\ldots,r\} $.The stationary distribution is 
\[\pi_{N,r}(j)=\frac{\binom{r}{j}\binom{N-r}{r-j}}{\binom{N}{r}}\]
and for $x\in {\mathbb X}$, the rates are given by
\[p_x=\frac{(r-x)(N-r-x)}{r(N-r)},\ \  q_x=\frac{x^2}{r(N-r)}\]
The eigenvalues of this chain are well known (see \cite{Karlin} or \cite {diacshas}) and are given by
\begin{equation}\label{eq:eigenberlap}
\beta_{N,r,i}=1-\frac{i(N-i+1)}{r(N-r)},\ 0\leq i\leq r
\end{equation}
The first eigenvectors are 
\[V_0\equiv 1,\ V_1(x)=(1-\frac{Nx}{r(N-r)})C_1\mbox{ where } C_1 \mbox{ is defined by the condition }\|V_1\|_{L^2(\pi_{N,r})}=1.\]

The sign of $V_1(0)V_1(x)$ is that of $1- \frac{Nx}{r(N-r)}$ so  the halting states belong to $\{x/ x\geq \frac{r(N-r)}{N}\}$.

In the symmetric case where $r=\frac{N}{2}$, a halting state must be $\geq\frac{r}{2}$.

\bprop\label{prop:Berlap}
$[\frac{N+1}{2}]$ is the smallest halting state for total variation.
\eprop

\begin{proof}
Unfortunatly, this chain is not monotone ($p_0+q_1>1$) so the proposition(\ref{prop:evenmonotone+}) cannot be applied.

Consider the case of $N$ even, and change $N$ to $2N$.

We use the notation of part (3.3.2) and of proposition (\ref{prop:evenmonotone+}). The result will be true if $Q_+^2$ is monotone and $x\mapsto \frac{Q_+^3(0,x)}{\widetilde{\pi}(x)}$ is decreasing.

Indeed, for all $n\geq 0$, then $x\mapsto \frac{(Q_+^3(0,.)Q_+^{2n})(x)}{\widetilde{\pi}(x)}=\frac{Q_+^{2n+3}(0,x)}{\widetilde{\pi}(x)}$ is decreasing and so $P^{2n+3}(0,N)=Q_+^{2n+3}(0,N)\leq \widetilde{\pi}(N)=\pi(N)$.

In the same way, using decreasing of $x\mapsto \frac{\delta_0}{\widetilde{\pi}}(x)$, for all $n\geq 0$, $P^{2n}(0,N)\leq \pi(N)$.

So first, we have to prove that for all $x$, for all $y$,
\begin{equation}\label{Q^2}
 Q_+^2(x,\{0,\ldots,y\})\geq Q_+^2(x+1,\{0,\ldots,y\})
\end{equation}

We have that for all $N-2\geq x\geq 1$, $p_x+q_{x+1}\leq 1$ and $p_{N-1}+2q_N\leq 1$, so if $x\geq 1$, $Q_+(x,\{0,\ldots,y\})\geq Q_+(x+1,\{0,\ldots,y\})$. This inequality is also true if $x=0$ and $y\geq 1$.

So if $x\geq 2$,
\[\begin{array}{ll}
 Q_+^2(x,\{0,\ldots,y\})&=q_x( Q_+(x-1,\{0,\ldots,y\})- Q_+(x,\{0,\ldots,y\}))\\
&\ \ +(1-p_x)(Q_+(x,\{0,\ldots,y\})-Q_+(x+1,\{0,\ldots,y\}))\\
&\ \ +Q_+(x+1,\{0,\ldots,y\})-Q_+(x+2,\{0,\ldots,y\})\\
&\ \ +Q_+(x+2,\{0,\ldots,y\})\\
&\geq q_{x+1}(Q_+(x,\{0,\ldots,y\})-Q_+(x+1,\{0,\ldots,y\}))\\
&\ \ +(1-p_{x+1})Q_+(x+1,\{0,\ldots,y\})-Q_+(x+2,\{0,\ldots,y\})\\
&\ \ +Q_+(x+2,\{0,\ldots,y\})=Q_+(x+1,\{0,\ldots,y\})
  \end{array}\]

It remains to prove (\ref{Q^2}) for $x=0$ and $x=1, y=0$.
\[\begin{array}{l}
 Q_+^2(0,\{0\})=p_0q_1+r_0^2=q_1\geq r_1q_1=Q_+^2(1,\{0\})\\
Q_+^2(0,\{0,1\})=q_1+r_1\\
Q_+^2(1,\{0,1\})=q_1+r_1(1-p_1)+p_1q_2=q_1+r_1+p_1(q_2-r_1)\\
Q_+^2(1,\{0\})=r_1q_1\\
Q_+^2(2,\{0\})=q_1q_2
\end{array}\]

So the condition is satisfied if $q_2\leq r_1$ and this the case for the Bernoulli-Laplace chain.

Now it remains to prove that $x\mapsto \frac{Q_+^3(0,x)}{\widetilde{\pi}(x)}$ is decreasing.

Denote $\widetilde{\pi}_n=Q_+^n(0,.)$.

We have by reversibility that $\frac{\widetilde{\pi}_{n+1}}{\widetilde{\pi}}(x)=\ds\sum_y\frac{\widetilde{\pi}_{n}}{\widetilde{\pi}}(y)Q_+(x,y)$.

So by some calculus, we find 
\[\frac{\widetilde{\pi}_3}{\widetilde{\pi}}(0)=\frac{r_1}{\widetilde{\pi}_1(1)}\ \frac{\widetilde{\pi}_3}{\widetilde{\pi}}(1)=\frac{q_1+r_1^2+p_1q_2}{\widetilde{\pi}_1(1)}\ \frac{\widetilde{\pi}_3}{\widetilde{\pi}}(2)=\frac{r_1q_2+q_2r_2}{\widetilde{\pi}_1(1)}\ \frac{\widetilde{\pi}_3}{\widetilde{\pi}}(3)=\frac{q_2q_3}{\widetilde{\pi}_1(1)}\]

So we have to verify if $r_1\geq q_1+r_1^2+p_1q_2\geq r_1q_2+q_2r_2\geq q_2q_3$. It is true if $N$ is big enough.

Consider now the case where $N$ is odd and replace $N$ by $2N+1$.

By the same ideas, it suffices to show that $Q_1^2$ and $P^2$ are monotone where $Q_1$ is given in proposition (\ref{prop:monotoneodd}). In fact $p_x+q_{x+1}\leq 1$ if $1\leq x\leq 2N-1$, but using symetry, $P^2(2N+1-x,\{0,\ldots,y\})=1-P^2(x,\{0,\ldots,2N-y\})$, so there is no additional verifications to do. We conclude as in proposition (\ref{prop:monotoneodd}).

\end{proof}
By the same method as used for the Ehrenfest chain, using (\ref{prop:Berlap}) and the knowledge of eigenvalues given by (\ref{eq:eigenberlap}), the symmetric Bernoulli-Laplace birth-and-death chain on $\{0,\ldots,N\}$ exhibits a cutoff at time $\frac{N\ln N}{4}$ with a window $N$.

\subsubsection{Asymmetric continuous-time Erhenfest process}
Let $Z_t=(Z^{(1)}_t,\ldots,Z^{(N)}_t)$  be a continuous-time Markov chain on $(\Z/2\Z)^N$ where $\left( (Z_t^{(i)})_{t\geq 0}\right) _{1\leq i\leq N}$ are iid continuous-time Markov chain on $\Z/2\Z$ with infinitesimal generator given by $\left(\begin{array}{cc} -\lambda&\lambda\\\mu&-\mu\end{array} \right)$. We have
\begin{equation}\label{semigroup}
 \begin{array}{l}
  \Pr_0(Z_t^1=0)=p_t(0,0)=\frac{\mu}{\lambda+\mu}+\frac{\lambda}{\lambda+\mu}e^{-t(\lambda+\mu)}\\
  \Pr_0(Z_t^1=1)=p_t(0,1)=\frac{\lambda}{\lambda+\mu}-\frac{\lambda}{\lambda+\mu}e^{-t(\lambda+\mu)}\\
  \Pr_1(Z_t^1=0)=p_t(1,0)=\frac{\mu}{\lambda+\mu}-\frac{\mu}{\lambda+\mu}e^{-t(\lambda+\mu)}\\
  \Pr_1(Z_t^1=1)=p_t(1,1)=\frac{\lambda}{\lambda+\mu}+\frac{\mu}{\lambda+\mu}e^{-t(\lambda+\mu)}\\
 \end{array}
\end{equation}

If  $X_t=|Z_t|=\displaystyle\sum_{i=1}^NZ_t^{(i)}$, then $(X_t)_{t\geq 0}$ is the asymmetric continuous-time  Ehrenfest Markov chain on $\{0,\ldots,N\}$. 

Furthermore if $0\leq i,j\leq N$, 
\begin{equation}\label{semigroupasymerhen}\Pr_i(X_t=j)=\displaystyle\sum_{k+l=j}\binom{i}{k}\binom{N-i}{l}p_t(1,1)^kp_t(1,0)^{i-k}p_t(0,1)^lp_t(0,0)^{N-i-l}
 \end{equation}

The stationnary distribution $\pi_{\lambda,\mu}$ is the binomial distribution ${\cal B}(N,\frac{\lambda}{\lambda+\mu})$.

So $\frac{\Pr(X_t=x)}{\pi(x)}=(1+\frac{\lambda}{\mu}e^{-t(\lambda+\mu)})^{N-x} (1-e^{-t(\lambda+\mu)})^{x}$.

So it is easy to prove that the minimal halting state  is $x^*_{\lambda,\mu}=\lceil\frac{\lambda}{\lambda+\mu}N\rceil$. 

\bprop
The asymmetric continuous-time Erhenfest process has a cutoff at time $\frac{\log(N)}{2(\lambda+\mu)}$.
\eprop
\begin{proof}

 The smallest halting state is given by $x^*_{\lambda,\mu}=\lceil\frac{\lambda}{\lambda+\mu}N\rceil$. Let $y^*=N-x^*$.
 
 We have to prove that the hypothesis of (\ref{prop:cutoffgeneral}) are satisfied.
 
 Denote by $Q_{\lambda,\mu}$ the generator of our chain and by $A$ the restriction of $Q_{\lambda,\mu}$ to $\{0,\ldots,x^*_{\lambda,\mu}-1\}$ and by $B$ the restriction of $Q_{\mu,\lambda}$ to $\{0,\ldots,y^*_{\mu,\lambda}-1\}$.
 
 For all $0\leq i\leq N$, $q_{\lambda,\mu}(i,i+1)=\lambda(N-i)$ and $q_{\lambda,\mu}(i,i-1)=\mu i$.
 
 Let $\varphi_0^{\lambda,\mu},\varphi_1^{\lambda,\mu},\ldots,\varphi_N^{\lambda,\mu}$ the sequence of orthogonal polynomials associated to $Q_{\lambda,\mu}$. 
 
 The eigenvalues of $-Q_{\lambda,\mu}$ are $t_k=k(\lambda+\mu)$, $0\leq k\leq N$.
 
 If  $L$ is  the $(N+1)\times(N+1)$ matrix defined by $L_{i,j}=1$ if $i+j=N$ and  $L_{i,j}=0$ if not, then $LQ_{\lambda,\mu}L=Q_{\mu,\lambda}$.
 
 So by an argument of symmetry, for all $0\leq k\leq N$, for all $0\leq x\leq N$, $\varphi_{N-x}^{\lambda,\mu}(t_k)=\varphi_N^{\lambda,\mu}(t_k)\varphi_x^{\lambda,\mu}(t_k)$. Furthermore, by (\ref{eq:signeigenvector}), $(-1)^k\varphi_N^{\lambda,\mu}(t_k)>0$.
 
 Let $V_A=0<\alpha_0<\alpha_1<\ldots<\alpha_{x^*-1}$ be the eigenvalues of $-A$.
 
 Let $V_B=0<\beta_0<\beta_1<\ldots<\beta_{y^*-1}$ be the eigenvalues of $-B$.
 
 \blemm
 For all $0\leq i\leq N-1$, the interval $[t_i,t_{i+1}]$ contains a single element of $V_A\cup V_B$. 
  \elemm
  \begin{proof}
  It is a consequence of Sturm's theorem : The number of sign changes in the sequence $\varphi_0^{\lambda,\mu}(t),\varphi_1^{\lambda,\mu}(t),\ldots,\varphi_{N}^{\lambda,\mu}(t)$ gives the number of eigenvalues $>t$. The convention is that if $\varphi_{k+1}(t)=0$ and $\varphi_k(t)\not=0$, the sign of $\varphi_{k+1}(t)$ is equal to the sign of $\varphi_k(t)$.
  
  Apply this to $t=-t_1$ : 
  
  The single change of sign can occur in the sequence $\varphi_0^{\lambda,\mu}(-t_1),\varphi_1^{\lambda,\mu}(-t_1),\ldots,\varphi_{x^*}^{\lambda,\mu}(-t_1)$, and in this case , $-t_1<-\alpha_0<0$ 
  
  or can occur in the sequence $(\varphi_{x^*}^{\lambda,\mu}(-t_1),\ldots,\varphi_{N}^{\lambda,\mu}(-t_1))=-(\varphi_{0}^{\mu,\lambda}(-t_1),\ldots,\varphi_{N-y^*}^{\mu,\lambda}(-t_1))$, and in this case $-t_1<-\beta_0<0$.
  
  This method is reiterated.
  \end{proof}
 But $\Es_0[T_{x^*}^{\lambda,\mu}]=\displaystyle\sum_{i=0}^{x^*-1}\frac{1}{\alpha_i}$ and $\Es_0[T_{y^*}^{\mu,\lambda}]=\displaystyle\sum_{i=0}^{y^*-1}\frac{1}{\beta_i}$. Denote the first mean by $U_N$ et the second by $V_N$.
 
 The lemma implies that $\displaystyle\sum_{i=1}^N\frac{1}{t_i}\leq U_N+V_N\leq \frac{1}{\alpha_0}+\frac{1}{\beta_0}+\displaystyle\sum_{i=1}^N\frac{1}{t_i}$.
 
 As in the discrete time, we have by (\ref{propertyeigenvalue}), $\alpha_0\geq \pi_{\lambda,\mu}(\{x^*+1,\cdots, N\})t_1$.
 
 But $\pi_{\lambda,\mu}(\{x^*+1,\cdots, N\})\sim \frac{1}{2}$, then $U_N+V_N\sim \displaystyle\sum_{i=1}^N\frac{1}{t_i}\sim \frac{log(N)}{\lambda+\mu}$.
 
 \blemm
 $(U_N-V_N)$ is bounded.
 \elemm
 
 \begin{proof}
 By using the  equality proved in \cite{AldousFill}, chapter 2, lemma 12,
 \[\Es_i[T_j]=\frac{Z_{j,j}-Z_{i,j}}{\pi(j)}, \ \mbox{with } Z_{i,j}=\int_0^{+\infty}(\Pr_i(X_t=j)-\pi(j))dt,\]
 we have 
 \[U_N-V_N=\frac{1}{\pi_{\lambda,\mu}(x^*)}\int_0^{+\infty}(\Pr_N(X_t^{\lambda,\mu}=x^*)-\Pr_0(X_t^{\lambda,\mu}=x^*))dt\]
 
 By using formula (\ref{semigroupasymerhen}), 
 \[U_N-V_N=\int_0^{+\infty}(1+\frac{\mu}{\lambda}e^{-t(\lambda+\mu)})^{x^*}(1-e^{-t(\lambda+\mu)})^{y^*}-(1-e^{-t(\lambda+\mu)})^{x^*}(1+\frac{\lambda}{\mu}e^{-t(\lambda+\mu)})^{y^*}\ dt\]
 
 After the change of variables $u=e^{-t(\lambda+\mu)}$, we obtain
 \[U_N-V_N=\frac{1}{\lambda+\mu}\int_0^1\frac{(1+\frac{\mu}{\lambda}u)^{x^*}(1-u)^{y^*}-(1-u)^{x^*}(1+\frac{\lambda}{\mu}u)^{y^*}}{u}\ du\]
 
 Some calculus gives the result.
 \end{proof}
 Then as  $U_N+V_N\sim \frac{log(N)}{\lambda+\mu}$ and $U_N-V_N$ is bounded, $U_N\sim\frac{log(N)}{2(\lambda+\mu)}$ and the result comes from proposition (\ref{prop:cutoffgeneral}).
  
 \end{proof}

\section{Appendix }

\subsection{ Birth-and-death chain and orthogonal polynomials}

A reference for orthogonal polynomial and birth and death chain can be found in Karlin and McGregor(\cite{KarlinMcGregor},\cite{Karlin}).

Denote by $P$ the transition matrix on  $\{0,\ldots,N\}$ and by $P_k$ the restriction of $P$ to $\{0,\ldots,k\}$ for $0\leq k\leq N$.

For $1\leq k\leq N$, denote by $\varphi_k$ the polynomial given by $\varphi_k(t)=\frac{1}{p_0\cdots p_{k-1}}\mbox{det}(tI_k-P_{k-1})$. Let $\varphi_0(t)=1$.

So $(\varphi_0,\ldots,\varphi_N)$ is a family of polynomials which satisfy the following recurrence equation :
\begin{equation}\label{eq:polyortho}
\varphi_k(t)=\frac{t-r_{k-1}}{p_{k-1}}\varphi_{k-1}(t)-\frac{q_{k-1}}{p_{k-1}}\varphi_{k-2}(t) 
\end{equation}

Let $1=\lambda_0>\lambda_1>\cdots>\lambda_N$ be the eigenvalues of $P$.

It exists $(\mu_0,\ldots,\mu_{N})\in\R_+^{N+1}$ such that  $(\varphi_0,\ldots,\varphi_N)$ are the orthogonal polynomials in $L^2([-1,1],\mu)$ where $\mu=\ds\sum_{k=0}^N\mu_k\delta_{\lambda_k}$.

Furthermore $\int_{-1}^1\varphi_k(t)^2\mu(dt)=\frac{\pi(0)}{\pi(k)}=\frac{q_1\cdots q_k}{p_0\cdots p_{k-1}}$, $\mu_0=\pi(0)$.

We have for all $n\geq 0$, for all $x,y\in {\mathbb X}$,
\begin{equation}\label{eq:spectraldec} 
P^n(x,y)=\frac{\pi(y)}{\pi(0)}\ds\sum_{k=0}^N\lambda_k^n\varphi_x(\lambda_k)\varphi_y(\lambda_k)\mu_k
\end{equation}

So the normalized in $L^2(\pi)$ eigenvector $V_k$ associated to the eigenvalue $\lambda_k$ is given by $V_k(x)=\sqrt{\frac{\mu_k}{\mu_0}}\varphi_x(\lambda_k)$.

The zeros of the orthogonal polynomials have the following interlacing property : if $m>n$, there is a zero of $\varphi_m$ between any two zeros of $\varphi_n$.

So there is a zero of $\varphi_N$ in each interval $]\lambda_{k+1},\lambda_k[$ and as $\varphi_N(1)=1$, we have :
\begin{equation}\label{eq:signeigenvector}
\forall 0\leq k\leq N,\ \ (-1)^k\varphi_N(\lambda_k)>0 \ \ (\mbox{or }(-1)^kV_k(N)>0))
\end{equation}

\paragraph*{The case of a symmetric birth-and-death chain on $\{0,\ldots,2N+1\}$ }

Denote by $Q_1$ the matrix equal to $P_N$ unless the entry $Q_1(N,N)$ that is $r_N+p_N$ and by $C_1(t)$ the  characteristic polynomial of $Q_1$.

Denote by $Q_2$ the matrix equal to $P_N$ unless the entry $Q_2(N,N)$ that is $r_N-p_N$ and by $C_2(t)$ the  characteristic polynomial of $Q_2$.

Let $1=\beta_0>\beta_1>\cdots>\beta_N$ be the eigenvalues of $Q_1$.

Let $\alpha_0>\alpha_1>\cdots>\alpha_N$ be the eigenvalues of $Q_2$.

Let $\eta_0>\eta_1>\cdots>\eta_{N-1}$ be the eigenvalues of $P_{N-1}$. 

Let $\gamma_0>\gamma_1>\cdots>\gamma_N$ be the eigenvalues of $P_{N}$.

So by the interlacing property, for all $0\leq k\leq N-1$, $\beta_{k+1}<\eta_k<\beta_k$, $\alpha_{k+1}<\eta_k<\alpha_k$ and $\gamma_{k+1}<\eta_k<\gamma_k$.

We have $(p_0\cdots p_{N-1})^{-1}C_1(t)=(t-r_N-p_N)\varphi_N(t)-q_N\varphi_{N-1}(t)$,   $(p_0\cdots p_{N-1})^{-1}C_2(t)=(t-r_N+p_N)\varphi_N(t)-q_N\varphi_{N-1}(t)$ and $\varphi_{N+1}(t)=(t-r_N)\varphi_N(t)-q_N\varphi_{N-1}(t)$.

So $(p_0\cdots p_{N-1})^{-1}C_1(\alpha_k)=-2p_N\varphi_N(\alpha_k)$. But as $\eta_k<\alpha_k<\eta_{k-1}$, $(-1)^k\varphi_N(\alpha_k)>0$.

But $\alpha_k\in ]\eta_k,\beta_k[$ or $\alpha_k\in ]\beta_k,\eta_{k-1}[$ and $(-1)^kC_1<0$ on the first interval and  $(-1)^kC_1>0$ on the second interval. So we have for all $0\leq k\leq N$, $\alpha_k<\beta_k$.

As $(p_0\cdots p_{N-1})^{-1}(C_1(t)+C_2(t))=2\varphi_{N+1}(t)$, $C_1(\gamma_k)=-C_2(\gamma_k)$ so for all $0\leq k\leq N-1$,

\begin{equation}\label{eq:eigenvaluesymm}
 \eta_k<\alpha_k<\gamma_k<\beta_k<\eta_{k-1}
\end{equation}

with $\eta_{-1}=+\infty$ et $\eta_N=-\infty$.

Suppose now that for all $x\in {\mathbb X}$, $r_x=0$.

In this case if $k$ is odd then $\varphi_k$ is odd and if $k$ is even, $\varphi_k$ is even.

Furthermore $C_2(-t)=(-1)^{N-1}C_1(t)$ and so for all $0\leq k\leq N$, $\alpha_k=-\beta_{N-k}$. So we have :
\begin{equation}\label{eq:r=0eigenvalue}
\left\lbrace \begin{array}{ll}\mbox{If }N=2a+1,& 0<|\beta_{a+1}|<\beta_a<|\beta_{a+2}|<\cdots<\beta_1<|\beta_{2a+1}|<\beta_0=1\\
\mbox{If }N=2a, &0<\beta_a<|\beta_{a+1}|<\beta_{a-1}<|\beta_{a+2}|<\cdots<\beta_1<|\beta_{2a}|<\beta_0=1\end{array}\right.
\end{equation}

So by (\ref{eq:spectraldec}), we have in the case where $N$ is even :
\begin{equation}\label{eq:r=0N+1}
P^{2n+1}(0,N+1)=\frac{2\pi(N+1)}{\pi(0)}\ds\sum_{k=0}^N\beta_k^{2n+1}\varphi_{N+1}(\beta_k)\mu_{2k}
\end{equation}
and $\varphi_{N+1}(\beta_k)$ has the sign of $(-1)^k$.

\subsection{Proof of results on continuous time}

For each $0\leq h$, consider the $h$-skeleton discrete time Markov chain $X^{(h)}(n)=X(nh)$, $n\in\N$. 

$X^{(h)}$ is a Markov chain with transition matrix $e^{hQ}$, initial distribution $\pi_0$ and stationary distribution $\pi$.

We denote the corresponding quantities (\ref{eq:psin}) and (\ref{Jn}) by $\psi_n^{(h)}$ and $J_n^{(h)}=\ds\prod_{k=0}^n(1-\psi_k^{(h)}(X_k^{(h)}))$

\bprop\label{prop:convJ}
For all $t\notin {\cal A}$, there exists $(h_n)\longrightarrow 0$ such that a.e. on $\Omega$
\[J^{(h_n)}_{[\frac{t}{h_n}]}\longrightarrow J_t\] 
\eprop

\begin{proof}

For all $x\in {\mathbb X}$, $t\longmapsto \pi_t(x)$ is analytic with derivative $\pi_tQ(x)$.

By consequence ${\cal A}_x=\{t\geq 0/ \pi_t(x)=\pi(x)\}$ is either equal to $\R_+$ or contains only isolated points.

Let ${\cal T}=\{x\in {\mathbb X}/\ \forall t\geq 0,\ \pi_t(x)=\pi(x)\}$ and  ${\cal A}=\ds\cup_{x\notin{\cal T}}{\cal A}_x$.

For all $t\geq 0$, ${\cal A}\cap [0,t]$ is a finite set.

Let $t\notin {\cal A}$. 
 We take $(h_n)\longrightarrow 0$ such that there exists $\epsilon>0$ such that for all $n\geq 0$, for all $a\in  {\cal A}\cap [0,t]$, $d(a,h_n\N)\leq h_n^{1+\epsilon}$.

We choose $\omega$ such that for all $r\geq 1$, $T_r(\omega)\notin {\cal A}\cup \{t\}$ where $(T_r)_{r\geq 0}$ are the jump times of the continuous-time Markov chain $(X_t)_{t\geq 0}$.

\begin{description}
 \item[First case ]: $\forall s\in[0,t], \pi_s(X_s(\omega))>\pi(X_s(\omega)$ 

So for all $s\in[0,t]$, $X_s(\omega)\notin {\cal T}$ and as $ T_r(\omega)\notin {\cal A}\cup \{t\}$, by continuity there exists $\alpha >0$ such that for all $s\in[0,t]$, $\pi_s(X_s(\omega))-\pi(X_s(\omega)>\alpha$.

\[ \ln J_n^{(h_n)}(\omega)=\ln \left(1-\frac{\pi(X_0(\omega))}{\pi_0(X_0(\omega))}\right)+\ds\sum_{k=1}^{[\frac{t}{h_n}]}\ln\left(1-\psi_k^{(h_n)}(X_{kh_n}(\omega) )\right)\]

So, we have, using continuity and derivability of $t\mapsto \pi_t(x)$ :
\[ \ln J_n^{(h_n)}(\omega)=\ln \left(1-\frac{\pi(X_0(\omega))}{\pi_0(X_0(\omega))}\right)+h_n\ds\sum_{k=1}^{[\frac{t}{h_n}]}\frac{(\pi_{h_nk}\wedge\pi)Q(X_{kh_n}(\omega) )}{\pi_{h_nk}(X_{kh_n}(\omega) )-\pi(X_{kh_n}(\omega) )}+O(h_n)\]

That gives the result in this case

\item[Second case ]: $\exists s\in[0,t], \pi_s(X_s(\omega))\leq\pi(X_s(\omega)$ 

If there exists $s\in [0,t]$ such that $\pi_s(X_s(\omega))<\pi(X_s(\omega))$ then by continuity, $J_n^{(h_n)}(\omega)=0$ for $n$ big enough.

So we suppose that there exists $s_0\in]0,t[$ which satisfies $\pi_{s_0}(X_{s_0}(\omega))=\pi(X_{s_0}(\omega))$ and for all $u\in [0,t]$, $\pi_{u}(X_{u}(\omega))\geq \pi(X_{u}(\omega))$.

Let $x=X_{s_0}(\omega)$ and $r$ such that $T_r(\omega)<s_0<T_{r+1}(\omega)$.

We want to prove that $J_n^{(h_n)}(\omega)\longrightarrow 0$.

Let $s_n\in h_n\N$ which satisfies $|s_n-s_0|\leq h_n^{1+\epsilon}$. For $n$ big enough, $T_r(\omega)<s_n-h_n<s_n<T_{r+1}(\omega)$.

We have $J_n^{(h_n)}(\omega)\leq \frac{\pi_{s_n}(x)-\pi(x)}{\pi_{s_n}(x)-(\pi_{s_n-h_n}\wedge\pi)e^{h_nQ}(x)}$.

Define for $n\geq 1$, 
\[R_n=\{x\in {\mathbb X}/ \pi_{s_0}(x)=\pi(x),\pi_{s_0}Q^k(x)=0\mbox{ for } 1\leq k\leq n-1,\ (\pi_{s_0}\wedge\pi)Q^k(x)=0\mbox{ for }  1\leq k\leq n\}\]

Let $n\geq 2$. We prove by iteration that 
\begin{equation}\label{R_n}
 x\in R_n\Longrightarrow \forall y/q(y,x)\not= 0,\ y\in R_{n-1},\ (\pi_{s_0}\wedge\pi)Q^n(y)\leq 0\mbox{ and }\pi_{s_0}Q^{n-1}(y)\geq 0
\end{equation}

Let $n\geq 2$. We prove by iteration that 
\begin{equation}\label{R_nbis}
 x\in R_n\mbox{ and } \pi_{s_0}Q^{n}(x)=0\Longrightarrow\forall y/q(y,x)\not= 0,y\in R_{n-1}\mbox{ and }\pi_{s_0}Q^{n-1}(y)=0
\end{equation}

For $n$ big enough, $\pi_{s_n}(x)-(\pi_{s_n-h_n}\wedge\pi)e^{h_nQ}(x)=\pi_{s_n}(x)-\pi(x)+\displaystyle\sum_{k\geq 1}(\pi_{s_n-h_n}\wedge\pi)Q^k(x)\frac{h_n^k}{k!}$.

Define the following sets :
\[{\cal P}=\{y\in {\mathbb X}/ \pi_{s_0}(y)>\pi(y)\},\ {\cal M}=\{y\in {\mathbb X}/ \pi_{s_0}(y)<\pi(y)\}\]
\[ {\cal A}_+=\{y\in {\mathbb X}/ \pi_{s_0}(y)=\pi(y), (\pi_{s_0}-\pi)(y-)>0\},\ {\cal A}_-=\{y\in {\mathbb X}/ \pi_{s_0}(y)=\pi(y), (\pi_{s_0}-\pi)(y-)<0\}\]
So, if $n$ is big enough and $s=s_n-h_n$,
\[(\pi_s\wedge\pi)Q^k(x)=\displaystyle\sum_{y\in {\cal T}\cup{\cal P}\cup{\cal A}_+}\pi(y)q_k(y,x)+\displaystyle\sum_{y\in {\cal M}\cup{\cal A}_-}\pi_s(y)q_k(y,x)=\displaystyle\sum_{y\in {\cal M}\cup{\cal A}_-}(\pi_s(y)-\pi(y))q_k(y,x)\].

So if we denote $\alpha_{l,k}=\displaystyle\sum_{y\in {\cal M}\cup{\cal A}_-}\pi_{s_0}Q^k(y)q_l(y,x)$, we have
\[(\pi_{s_n-h_n}\wedge\pi)Q^{l}(x)=(\pi_{s_0}\wedge\pi)Q^{l}(x)+\displaystyle\sum_{r\geq 1}\frac{(s_n-h_n-s_0)^r}{r!}\alpha_{l,r}\ \ \mbox{and so}\]

\begin{equation}\label{dl}
\pi(x)-(\pi_{s_n-h_n}\wedge\pi)e^{h_nQ}(x)=-\displaystyle\sum_{l\geq 1}(\pi_{s_0}\wedge\pi)Q^{l}(x)\frac{h_n^l}{l!}-\displaystyle\sum_{l\geq 1}\frac{h_n^l}{l!}\displaystyle\sum_{r\geq 1}\frac{(s_n-h_n-s_0)^r}{r!}\alpha_{l,r}
\end{equation}

We prove by induction that for $n\geq 2$,
\begin{equation}\label{alpha1} 
x\in R_n\Longrightarrow\left\lbrace\begin{array}{ll}
                                  \mbox{if }  k+l\leq n-1 & \mbox{then }\alpha_{k,l}=0\\
                                  \mbox{if }  k+l= n & \mbox{then }\alpha_{k,l}\geq 0
                                   \end{array}
 \right.
\end{equation}

We prove  by induction that for $n\geq 2$, 
\begin{equation}\label{alpha2} 
x\in R_n\mbox{ and } \pi_{s_0}Q^n(x)=0\Longrightarrow\left\lbrace\begin{array}{ll}
                                  \mbox{if }  k+l\leq n & \mbox{then }\alpha_{k,l}=0\\
                                  \mbox{if }  k+l= n+1 & \mbox{then }(-1)^l\alpha_{k,l}\leq 0
                                   \end{array}
 \right.
\end{equation}

Let $2k_0$ be the smallest integer $k$ which satisfies $\pi_{s_0}Q^k(x)\not=0$. We have $\pi_{s_0}Q^{2k_0}(x)>0$.

Let $k_1$ be the smallest integer $k$ which satisfies $(\pi_{s_0}\wedge\pi)Q^k(x)\not=0$.

If $k_1\leq 2k_0$ then by (\ref{dl})
\[\pi(x)-(\pi_{s_n-h_n}\wedge\pi)e^{h_nQ}(x)=-(\pi_{s_0}\wedge\pi)Q^{k_1}(x)\frac{h_n^{k_1}}{{k_1}!}-\displaystyle\sum_{l,r\geq 1, l+r\leq k_1}\frac{h_n^l}{l!}\frac{(s_n-h_n-s_0)^r}{r!}\alpha_{l,r}+o(h_n^{k_1})\]

But as $x\in R_{k_1-1}$, using (\ref{R_n}), we have $(\pi_{s_0}\wedge\pi)Q^{k-1}(x)<0$ and using (\ref{alpha2}) $(-1)^r\alpha_{l,r}(x)\leq 0$ if $l+r=k_1$ and $\alpha_{l,r}(x)=0$ if $l+r<k_1$, so there exists $C>0$ such that 
\[\pi(x)-(\pi_{s_n-h_n}\wedge\pi)e^{h_nQ}(x)=h_n^{k_1}C+o(h_n^{k_1})\]

If $k_1>2k_0$, then $x\in R_{2k_0}$. By (\ref{dl}), we have
\[\pi(x)-(\pi_{s_n-h_n}\wedge\pi)e^{h_nQ}(x)=h_n^{2k_0}\left(\displaystyle\sum_{l,r\geq 1, l+r=2k_0}(-1)^{r+1}\frac{\alpha_{l,r}}{l!r!} \right)+o(h_n^{2k_0})\]

But $\alpha_{1,2k_0-1}(x)=\displaystyle\sum_{y\in{\cal A}_-}(\pi_{s_0}Q^{2k_0-1})(y)q(y,x)$. By (\ref{R_n}) and (\ref{R_nbis}), if $q(y,x)>0$, then $y\in R_{2k_0-2}$ and $\pi_{s_0}Q^{2k_0-1}(y)\geq 0$.

Furthermore if $y\in {\cal A}_+$, then  $\pi_{s_0}Q^{2k_0-1}(y)= 0$. So $\alpha_{1,2k_0-1}(x)=\pi_{s_0}Q^{2k_0}(x)>0$.

And using again (\ref{alpha2}), there exists $C>0$ such that 

\[\pi(x)-(\pi_{s_n-h_n}\wedge\pi)e^{h_nQ}(x)=h_n^{2k_0}C+o(h_n^{2k_0})\]

The fact that $\pi_{s_n}(x)-\pi(x)=\pi_{s_0}Q^{2k_0}(x)\frac{(s_n-s_0)^{2k_0}}{(2k_0)!}+o((s_n-s_0)^{2k_0})$ and $|s_n-s_0|\leq h_n^{1+\epsilon}$ implies that $\frac{\pi_{s_n}(x)-\pi(x)}{\pi_{s_n}(x)-(\pi_{s_n-h_n}\wedge\pi)e^{h_nQ}(x)}$ goes to $0$.
\end{description}
\end{proof}

We give now the proof of proposition (\ref{prop:weakcontcoup}).

\begin{proof}

We define for $t\geq 0$ and $x \in {\mathbb X}$, $f(t,x)= \Pr(T\leq t, X_T=x)$ and $l(t,x)=\Pr(T\leq t, X_t=x)$.
This proposition is the consequence of the two following lemma.

 \blemm
$l(t,x)=(\pi_t\wedge\pi)(x)$
\elemm
\begin{proof}
By using the proposition (\ref{prop:convJ}), we have 
\[l(t,x)=\Pr(U\geq J_t, X_t=x)=\ds\lim_{n\rightarrow +\infty}\Pr(U\geq J^{(h_n)}_{[\frac{t}{h_n}]},X_t=x)\]

But $h_n [\frac{t}{h_n}]\leq t$ so, 
\[\begin{array}{ll}\Pr(U\geq J^{(h_n)}_{[\frac{t}{h_n}]},X_t=x)&=\sum_{y\in {\mathbb X}}\Pr(U\geq J^{(h_n)}_{[\frac{t}{h_n}]},X_{h_n [\frac{t}{h_n}]}=y,X_t=x)\\
 &=\sum_{y\in {\mathbb X}}\Pr(U\geq J^{(h_n)}_{[\frac{t}{h_n}]},X_{h_n [\frac{t}{h_n}]}=y)e^{(t-h_n[\frac{t}{h_n}])Q}(y,x)\\
&=\sum_{y\in {\mathbb X}}(\pi_{h_n [\frac{t}{h_n}]}\wedge\pi)(y)e^{(t-h_n[\frac{t}{h_n}])Q}(y,x)\\
&\rightarrow (\pi_t\wedge\pi)(x)\end{array}\]

this by (\ref{loiXTT}) and (\ref{loiXTT2}).
\end{proof}

\blemm
$\frac{\partial f}{\partial t}(t,x)=\frac{\partial l}{\partial t}(t,x)-\ds\sum_{y\in {\mathbb X}}l(t,y)q(y,x)$
\elemm

It is true for all stopping times.

\end{proof}
\subsection{Proof of lemma (\ref{le:le1})}
\blemm
Let $\Gamma_n=\ds\sum_{k=0}^{2d}\alpha^n_kB_k$ where
\begin{equation}\label{le:conGamma}
B_0>0,\ \forall i B_{2i}B_{2i+1}<0,\ B_{2i+1}B_{2i+2}>0,\ \alpha_0= 1>\alpha_1>\alpha_2>\ldots>\alpha_{2d}>0
\end{equation}
We suppose that $\Gamma_0=\Gamma_1=\cdots=\Gamma_{d-1}=0$. Then $n\mapsto \Gamma_n$ is stricly increasing on $\{d,d+1,\cdots\}$.

\elemm

\begin{proof}
We prove the result inductively by the study of  the real function $\varphi(t)=B_0+\ds\sum_{k=1}^{2d}B_ke^{-\mu_kt}$ where $\mu_k=-\ln\alpha_k$.

We prove that if (\ref{le:conGamma}) is satisfied then \begin{itemize}
                                                     \item 
                                                    
If $d$ is even, there exists $1\leq k\leq \frac{d}{2}$, there exists $t_1<\ldots<t_{2k-1}$ such that

$\forall 1\leq i\leq k$, $\varphi$ is strictly increasing on $[t_{2i-1},t_{2i}[$

$\forall 0\leq i\leq k-1$, $\varphi$ is strictly decreasing on $[t_{2i},t_{2i+1}[$

where $t_0=-\infty, t_{2k}=+\infty$.

\item If $d$ is odd, there exists $1\leq k\leq \frac{d-1}{2}$, there exists $t_1<\ldots<t_{2k}$ such that

$\forall 0\leq i\leq k$, $\varphi$ is strictly increasing on $[t_{2i},t_{2i+1}[$

$\forall 1\leq i\leq k$, $\varphi$ is strictly decreasing on $[t_{2i-1},t_{2i}[$

where $t_0=-\infty, t_{2k+1}=+\infty$.
\end{itemize}

It is easy to see that if $d=1$, $\varphi$ is strictly increasing.

If $d$ is odd and $\varphi(t)=B_0+\ds\sum_{k=1}^{2(d+1)}B_ke^{-\mu_kt}$, we write $\varphi'(t)=-e^{-\mu_1t}\psi(t)$ and $\psi'(t)=e^{-(\mu_2-\mu_1)t}\theta(t)$.

$\theta $ satisfy the same hypothesis so there exists $1\leq k\leq \frac{d-1}{2}$, there exists $t_1<\cdots<t_{2k}$ such that

$\forall 0\leq i\leq k$, $\theta$ is strictly increasing on $[t_{2i},t_{2i+1}[$

$\forall 1\leq i\leq k$, $\theta$ is strictly decreasing on $[t_{2i-1},t_{2i}[$

where $t_0=-\infty, t_{2k+1}=+\infty$.

Furthermore $\ds\lim_{t\rightarrow +\infty}\theta(t)=-B_2\mu_2(\mu_2-\mu_1)>0$.

The number of zeros of $\theta$ is odd and smaller than $2k+1$ . So there exists $l\leq k$, there exists $s_1<\cdots <s_{2l+1}$ such that 

$\theta<0$ on $]s_{2i},s_{2i+1}[$ for $0\leq i\leq l$, $s_0=-\infty$.

$\theta>0$ on $]s_{2i+1},s_{2i+2}[$ for $0\leq i\leq l$, $s_{2l+2}=+\infty$.

Furthermore $\ds\lim_{t\rightarrow +\infty}\psi(t)=B_1\mu_1<0$.

So the number of zeros of $\psi$ is odd and smaller than $2l+1$. So there exists $m\leq l$, there exists $r_1<\cdots <r_{2m+1}$ such that 

$\psi>0$ on $]r_{2i},r_{2i+1}[$ for $0\leq i\leq m$, $r_0=-\infty$.

$\psi<0$ on $]r_{2i+1},r_{2i+2}[$ for $0\leq i\leq m$, $r_{2m+2}=+\infty$.

As $m\leq l\leq k\leq\frac{d-1}{2}$, we have proved the result for $\varphi$.

The case where $d$ is even is the same. 

The preceding result prove that $\varphi(t)=B_0+\ds\sum_{k=1}^{2d}B_ke^{-\mu_kt}$ has at most $d$ zeros. So if  $\varphi(0)=\varphi(1)=\cdots=\varphi(d-1)=0$, $t\mapsto \varphi(t)$ is strictly increasing on $]d-1,+\infty[$ and this proves the lemma.

\end{proof}

\subsection{Proof of lemma (\ref{le:inc})}

\blemm 
Let $1=\lambda_0>\lambda_1>\ldots>\lambda_d>0$ and $\Gamma_n=\ds\sum_{i=0}^dA_i\lambda_i^n,A_0\geq 0$

If for $0\leq i\leq d-1$, $\Gamma_i=0$, then $n\rightarrow \Gamma_n$ is increasing and so for all $n\geq 0$, $\Gamma_n\leq A_0$.
\elemm 

\begin{proof}
We show the result inductively.

If $d=1$, $\Gamma_n=A_0(1-\lambda_1^n)$.

We suppose the result true for $d-1$. We can suppose that for all $0\leq i\leq d$, $A_i\not= 0$ et denote by $\epsilon$ the sign of $A_1$.

Then for all $n\geq 0$, $\Gamma_{n+1}-\Gamma_n=\epsilon\lambda_1^n\theta_n$, where
\[\theta_n=B_0+\sum_{i=1}^{d-1}B_i\widetilde{\lambda}_i^n,\ \ B_i=\epsilon A_{i+1}(\lambda_{i+1}-1),\ \widetilde{\lambda}_i=\frac{\lambda_{i+1}}{\lambda_1}\]
Furthermore, 
$\theta_0=\cdots=\theta_{d-2}=0$, so by hypothesis, $n\rightarrow \theta_n$ is increasing, so as $\theta_0=0$, for all $n\geq 0$, $\theta_n\geq 0$ and $n\rightarrow \epsilon\Gamma_n$ is increasing.

So $0=\epsilon\Gamma_0\leq \epsilon\Gamma_n\leq \lim_{n\rightarrow +\infty}\epsilon\Gamma_n=\epsilon A_0$.

As $A_0>0$, then $\epsilon=1$. this completes the proof of the lemma.
\end{proof}

% \bibliographystyle{plain}
% \bibliography{bibliomax}

\end{document}